\documentclass[journal]{IEEEtran}
\usepackage{amsmath,amsfonts}
\usepackage{algorithmic}
\usepackage{algorithm}
\usepackage{array}
\usepackage[caption=false,font=normalsize,labelfont=sf,textfont=sf]{subfig}
\usepackage{textcomp}
\usepackage{stfloats}
\usepackage{url}
\usepackage{verbatim}
\usepackage{graphicx}
\usepackage{cite}
\usepackage{xcolor}
\usepackage{tabularx}
\newcolumntype{Y}{>{\centering\arraybackslash}X}

\begin{document}

\title{Customized Interior-Point Methods Solver\\ for Embedded Real-Time Convex Optimization}

\author{Jae-Il Jang and Chang-Hun Lee\textsuperscript{*}
        % <-this % stops a space
\thanks{This paper is based on work contained in the author's Master's thesis submitted to Korea Advanced Institute of Science and Technology (KAIST) in February 2025. (\textit{Corresponding author : Chang-Hun Lee}.)}
\thanks{Jae-Il Jang, and Chang-Hun Lee are with the Department of Aerospace Engineering, Korea Advanced Institute of Science and Technology (KAIST), 291, Daehak-ro, Yuseong-gu, Daejeon, Republic of Korea. (E-mail: jji249@kaist.ac.kr, lckdgns@kaist.ac.kr). }
}

\markboth{}%
{Shell \MakeLowercase{\textit{et al.}}: A Sample Article Using IEEEtran.cls for IEEE Journals}

%\IEEEpubid{0000--0000/00\$00.00~\copyright~2025 IEEE}

\maketitle

\begin{abstract}
    This paper presents a customized second-order cone programming (SOCP) solver tailored for embedded real-time optimization, which frequently arises in modern guidance and control (G\&C) applications.
    The solver employs a practically efficient predictor-corrector type primal-dual interior-point method (PDIPM) combined with a homogeneous embedding framework for infeasibility detection. Unlike conventional homogeneous self-dual embedding formulations, the adopted approach can directly handle quadratic cost functions without requiring problem reformulation. 
    This capability allows the solver to directly address quadratic objective SOCP problems, while avoiding unnecessary performance degradation caused by the loss of sparsity due to problem reformulation.
    To support a systematic workflow, we also develop a code generation tool that analyzes the sparsity pattern of the problem to be solved and generates customized solver code using a predefined code template. The generated solver code is written in C with no external dependencies other than the standard library \texttt{math.h}, and it supports complete static allocation of all data.
    Additionally, it provides parsing information to facilitate the use of the solver by end users.
    Finally, benchmark and numerical experiments on an embedded platform demonstrate that the developed solver outperforms the existing solvers on problem scales typical of G\&C applications.
\end{abstract}


\section{Introduction}

    \IEEEPARstart{O}{ptimization} has become an important tool in modern guidance and control (G\&C) applications \cite{ferreau2017embedded}. By transcribing the G\&C problem into a constrained optimization problem, one can efficiently obtain desired solutions while accounting for complex dynamics, actuator limitations, and operational constraints. Among various approaches, convex optimization has emerged as a prominent approach due to its favorable theoretical and practical properties for onboard real-time applications, such as guaranteed convergence to a global optimum and existence of polynomial-time algorithms \cite{boyd2004convex, wang2024survey}. As a result, it plays a key role in enabling practical implementation of optimization-based techniques, such as model predictive control (MPC) \cite{wang2010fast, di2018real, di2018dynamic, hong2019model} and onboard trajectory optimization \cite{malyuta2022convex, acikmese2007convex, scharf2017implementation, blackmore2016autonomous, chen2023fast, wang2017constrained, liu2016exact, wang2018minimum}. It also frequently arises as subproblems in sequential methods for nonconvex problems, such as sequential quadratic programming (SQP) \cite{boggs1995sequential} and sequential convex programming (SCP) \cite{mao2016successive, bonalli2023analysis}.

    However, solving optimization problems onboard remains a major challenge to the practical implementation of these methods due to the stringent real-time execution requirements under embedded hardware with limited computational and memory resources \cite{ferreau2017embedded, dueri2017customized, chen2023fast}.
    Therefore, an efficient and reliable solver is essential for the successful deployment. Several numerical solution techniques have been developed over the past decades. In particular, interior-point methods (IPMs) \cite{nesterov1994interior,wright1997primal} and alternating direction method of multipliers (ADMM) \cite{boyd2011distributed} are among the most popular choices. Both algorithms can effectively handle general convex conic programs, including linear programs (LPs), quadratic programs (QPs), and second-order cone programs (SOCPs), which encompass most G\&C problems \cite{malyuta2021advances}. ADMM is particularly attractive for large-scale problems due to its low computational cost per iteration. However, it is known to suffer from poor tail convergence at high precision \cite{he2012on, boyd2011distributed}. It can also be sensitive to the ill-conditioned problem data, potentially requiring careful tuning of the parameters, such as step size, prior to execution \cite{themelis2019supermann}. These limitations may hinder reliable real-time embedded operation, and may adversely affect the convergence of sequential methods such as SQP and SCP \cite{leibfritz1999inexact}.

    \IEEEpubidadjcol
    
    On the other hand, IPMs are capable of achieving high-precision solutions in far fewer iterations, almost independent of problem size or conditioning. Moreover, since most G\&C problems are typically formulated as small- to medium-scale problems, IPMs tend to exhibit strong performance in such applications \cite[\S1]{ryu2022large}. They have been implemented in several modern convex optimization solvers, including CVXGEN \cite{mattingley2012cvxgen}, QOCO \cite{govind2025qoco}, ECOS \cite{domahidi2013ECOS}, MOSEK \cite{mosek2024}, and Clarabel \cite{goulart2024clarabel}. However, since CVXGEN and QOCO rely solely on standard IPM implementations, which assume the existence of a strictly feasible solution, their algorithms fail to converge when solving infeasible problem instances. Such infeasibility scenarios are commonly encountered in practical G\&C applications. For example, instances of the powered descent guidance problem may become infeasible depending on the specified time-of-flight \cite{acikmese2007convex, dueri2017customized}. To ensure reliable G\&C operation, it is important to detect infeasibility in problem instances within a finite number of iterations. A widely adopted approach to address this challenge is homogeneous embedding \cite{ye1994An, andersen1999homogeneous, yoshise2007interior}, which reformulates the problem into a slightly larger feasibility problem. This lifted problem is always feasible, which ensures that the IPM algorithm converges within a finite number of iterations regardless of the feasibility of the original problem. Once a solution is found from the lifted problem, it can be used either to recover the optimal solution to the original problem or to provide a certificate of infeasibility, depending on the status of the solution. In the case of ECOS and MOSEK, a homogeneous self-dual embedding \cite{ye1994An, vandenberghe2010cvxopt} formulation is employed to enable infeasibility detection. However, this formulation is limited to problems with linear cost functions due to the reliance on the self-dual formulation. Thus, to apply such a homogeneous embedding to problems with quadratic cost functions, the original problem needs to be reformulated as a linear objective problem by introducing a slack variable and an additional second-order cone constraint \cite{goulart2024clarabel}. This transformation leads to unnecessary increases in problem size and degrades the sparsity of the problem data, which can adversely affect solver performance. To address this issue, previous studies \cite{andersen1999homogeneous, odonoghue2021operator, yoshise2007interior} have proposed a homogeneous embedding formulation capable of handling quadratic cost functions. Clarabel is one of the few solvers that adopt this technique within an IPM framework and has demonstrated its effectiveness across diverse problem classes. Furthermore, it supports more general conic constraints beyond SOCPs, including semidefinite cones, exponential cones and power cones. However, it is designed as a general-purpose solver, rather than tailored for embedded applications. In particular, its implementation in Julia and Rust poses practical challenges for deployment on embedded systems. Although Rust is capable of being used in embedded environments and ensures memory safety through its ownership model, its ecosystem is still far less mature than that of C. As a result, limitations exist in terms of compatibility with existing C-based legacy systems and the availability of library support. 
    
    Another important factor in embedded real-time optimization is effectively exploiting the sparsity of problem data. Embedded real-time applications typically involve repeatedly solving problems with the same structure and sparsity pattern but with different data instances. This feature makes a \textit{custom solver} particularly well-suited for such applications. 
    A custom solver refers to a solver specifically tailored to a set of problems that share a common structure and sparsity pattern.
    It exploits the known structure and sparsity pattern to precompute certain information, such as the sparsity information for matrix factorization, during an offline phase. The precomputed information is then utilized during online operation, thereby improving computational efficiency. 
    Representative customized IPM solvers such as CVXGEN and QOCO employ custom linear algebra routines in which matrix operations and factorizations are explicitly hard-coded. While this approach enables faster operation by reducing overhead during execution, it has the drawback that the resulting code size increases rapidly with the problem size, limiting its scalability beyond small-scale problems \cite{dueri2014automated}. This issue can be mitigated by implementing the numerical routines for linear algebra using a loop-based structure, while caching only the sparsity information obtained offline \cite{domahidi2012efficient, banjac2017embedded}. A further distinctive feature of a customized solver is its ability to statically allocate all data. In the case of general-purpose solvers, dynamic memory allocation is inevitable in order to accommodate arbitrary problem structures. However, dynamic memory allocation introduces risks such as memory leaks and fragmentation, which are particularly undesirable in embedded systems where computational and memory resources are typically limited. Furthermore, it complicates the verification and validation process due to uncertainties introduced during the allocation process. In contrast, a customized solver can statically allocate all problem data by leveraging the known problem structure in advance, making it well-suited for embedded applications.

    Motivated by these observations, this study aims to develop a solver suitable for real-time embedded optimization. 
    To achieve efficiency and robustness, the solver developed in this study adopts the predictor-corrector type primal-dual IPM (PDIPM) \cite{mehrotra1992on} as its core algorithm, which is known to be practically efficient. Furthermore, the homogeneous embedding framework proposed in \cite{yoshise2007interior} is incorporated to directly handle SOCPs with quadratic cost functions while enabling infeasibility detection.
    To implement this algorithm in a form suitable for operations on embedded platforms, we employed a custom solver framework and developed a dedicated code generation tool. The problem structure and sparsity pattern are analyzed offline to improve computational efficiency during online execution. A customized solver written in C is automatically generated using the code generation tool, incorporating precomputed sparsity information. 
    
    The main contributions of this paper include the following:
    \begin{itemize}
        \item The solver developed in this study employs a PDIPM capable of infeasibility detection for a broad class of SOCPs, including LPs, QPs, SOCPs with linear objectives and even those with quadratic objectives. 
        This algorithmic formulation allows efficient treatment of a wide range of problems without reformulation, thereby avoiding unnecessary performance degradation.

        \item A solver customization framework is applied to enable efficient execution in embedded system environments. Moreover, a dedicated code generation tool is developed to reduce user workload by automating the customization process. The generated code has no external dependencies other than the standard library \texttt{math.h}, and its size remains nearly constant regardless of problem dimensions due to the use of loop-based structures for linear algebra routines. It also supports complete static allocation of problem data, and is designed for memory efficiency by utilizing only the upper triangular part of a symmetric matrix. In addition, the code generation tool provides parsing information to facilitate the use of the solver, significantly reducing the user workload for handling problem data in C-based environments.

        \item Extensive performance evaluations are conducted to demonstrate the effectiveness of the developed solver. Its performance and robustness are compared against existing solvers across standard benchmark problem sets. Furthermore, numerical experiments on an embedded system using G\&C examples validate its capability for real-time optimization in resource-constrained environments.
    \end{itemize}

    The remainder of this paper is organized as follows: Section \ref{sec: SOCP with quadratic cost} presents the problem setup addressed by the solver developed in this study. Section \ref{sec: solver algorithm} details the core algorithms of the solver, and Section \ref{sec: solver customization} describes the customization and code generation strategies. Section \ref{sec: benchmark test} shows the results of benchmark tests and Section \ref{sec: embedded system experiments} provides experiment results conducted on the embedded system with ARM Cortex-A9 processor. Finally, Section \ref{sec: conclusion} summarizes our conclusions.

\section{Second-Order Cone Programming with a Quadratic Cost Function} \label{sec: SOCP with quadratic cost}

    The customized solver developed in this study addresses the SOCP problem with a quadratic cost function expressed in the following form:
    \begin{equation}
        \begin{split}
            \underset{x, s}{\text{minimize}} & \quad \frac{1}{2}x^{T}Qx + q^{T}x \\
            \text{subject to} & \quad Ax = b \\
            & \quad Gx + s = h \\
            & \quad s \in \mathcal{K}
        \end{split}
        \tag{$\mathcal{P}$}
        \label{equ2: primal problem}
    \end{equation}
    \noindent
    where $x \in \mathbb{R}^{n}$ is the primal decision variable, $s \in \mathbb{R}^{m}$ is the slack variable, and $Q \in \mathbb{R}^{n \times n}$, $q \in \mathbb{R}^{n}$, $A \in \mathbb{R}^{p \times n}$, $b \in \mathbb{R}^{p}$, $G \in \mathbb{R}^{m\times n}$, $h \in \mathbb{R}^{m}$ represent the problem data. It is assumed that $Q$ is symmetric positive semidefinite, possibly zero, and that $A$ and $\begin{bmatrix} A^{T} & G^{T} \end{bmatrix}^{T}$ have full row rank. The cone $\mathcal{K}$ is a closed convex cone, and since the solver developed in this paper is designed to handle problems within the SOCP category, it is assumed that $\mathcal{K}$ consists of the nonnegative cone (NNC), 
    \begin{equation}
        \mathbb{R}_{+} := \left\{ x \in \mathbb{R} \left. \right\vert \ x \geq 0 \right\}
        \label{equ2: def nnc}
    \end{equation}
    \noindent
    the second-order cone (SOC) with dimension $d > 1$,
    \begin{equation}
        \mathbb{Q}^{d} :=  \left\{ (x_{0}, x_{1} ) \in \mathbb{R} \times \mathbb{R}^{d-1} \left.\right \vert \ \Vert x_{1} \Vert_{2} \leq x_{0}  \right\}
        \label{equ2: def soc}
    \end{equation}
    \noindent
    or both. Without loss of generality, it is further assumed that the cone $\mathcal{K}$ is ordered such that the first $l\geq 0$ cones in $\mathcal{K}$ correspond to the NNC, followed by $n_{SOC} \geq 0$ SOCs with dimension $d_{i}$ $(i=1, 2, \cdots, n_{SOC})$, which is given by
    \begin{equation}
        \mathcal{K} = \mathbb{R}_{+, 1} \times \cdots \times \mathbb{R}_{+, l} \times \mathbb{Q}^{d_{1}} \times \cdots \times \mathbb{Q}^{d_{n_{SOC}}}
        \label{equ2 : cone K}
    \end{equation}
    In the case where $l=0$ or $n_{SOC} = 0$, the corresponding constituent cones (i.e., the NNC blocks or the SOC blocks) are omitted in the representation of $\mathcal{K}$ in \eqref{equ2 : cone K}.
    \noindent
    The optimal solution to the primal problem \ref{equ2: primal problem}, if it exists, is denoted as $(x^{\ast}, s^{\ast})$, and the corresponding finite optimal value as $p^{\ast}$.
    
    The associated Lagrange dual problem to \ref{equ2: primal problem} can be written as follows:
    \begin{equation}
        \begin{split}
            \underset{x, y, z}{\text{maximize}} & \quad -\frac{1}{2}x^{T}Qx - b^{T}y - h^{T}z \\
            \text{subject to}&\quad Qx + A^{T}y + G^{T}z + q = 0 \\
            & \quad z \in \mathcal{K}^{\ast}
        \end{split}
        \tag{$\mathcal{D}$}
        \label{equ2: dual problem}
    \end{equation}
    \noindent
    where $y \in \mathbb{R}^{p}$, $z \in \mathbb{R}^{m}$ are dual variables, and $\mathcal{K}^{\ast}$ is the dual cone of $\mathcal{K}$, defined by
    \begin{equation}
        \mathcal{K}^{\ast} = \{ y \left. \right\vert \ \langle y, x \rangle \geq 0, \ \forall x \in \mathcal{K} \} 
    \end{equation}
    The optimal solution to the dual problem \ref{equ2: dual problem}, if it exists, is denoted as $(y^{\ast}, z^{\ast})$, with the corresponding finite optimal value $d^{\ast}$.
    
    Assume that the strong duality holds between problem \ref{equ2: primal problem} and \ref{equ2: dual problem}, i.e., $p^{\ast} = d^{\ast}$, which can be achieved by satisfying the mild regularity condition called Slater's condition \cite{boyd2004convex}. Then, the Karush-Kuhn-Tucker (KKT) conditions provide the necessary and sufficient conditions for optimality in a convex optimization problem \cite{boyd2004convex}. The KKT conditions for the given problem \ref{equ2: primal problem} and its dual \ref{equ2: dual problem} can be written as:
    \begin{subequations}
        \begin{align}
            Qx + A^{T}y + G^{T}z + q &= 0   \label{equ2: KKT conditions -a}  \\
            Ax &= b                         \label{equ2: KKT conditions -b}  \\
            Gx + s &= h                     \label{equ2: KKT conditions -c}  \\
            s^{T}z &= 0                     \label{equ2: KKT conditions -d}  \\
            (s, z) &\in \mathcal{K} \times \mathcal{K}^{\ast} \label{equ2: KKT conditions -e}
        \end{align}
        \label{equ2: KKT conditions}
    \end{subequations}
    \noindent
    Therefore, obtaining the primal-dual optimum $(x^{\ast}, s^{\ast}, y^{\ast}, z^{\ast})$ reduces to finding a primal-dual feasible $(x, s, y, z)$ with zero duality gap, since $s^{T}z$ in \eqref{equ2: KKT conditions -d} is identical to the duality gap for primal-dual feasible $(x, s, y, z)$ \cite[\S5.5.3]{boyd2004convex}.
    
\section{Solver Algorithm} \label{sec: solver algorithm}

\subsection{Interior-Point Methods}
    IPMs are among the most efficient algorithms for solving convex optimization problems with inequality constraints. They introduce a barrier function to incorporate the inequality constraints into the cost function. This modification transforms the original inequality-constrained problem into an equality-constrained one. They then iteratively solve this equality-constrained problem using Newton's method and converge to the optimal solution. A typical barrier function for IPMs over the given cone $\mathcal{K}$ is a log-det barrier function, which is defined as:
    \begin{equation}
        \Phi(x) = - \sum_{i = 1}^{l+n_{SOC}} \log \det (x_{i})
        \label{equ3: log-det barrier}
    \end{equation}
    where $x_{i}$ belongs to either the NNC or the SOC.
    \begin{equation}
        \det{(x)} = \begin{cases} x & x \in \mathbb{R}_{+} \\ x_{0}^{2} - x_{1}^{T}x_{1} & x \in \mathbb{Q}^{d} \end{cases}
    \end{equation}
    %Here, $N = l + n_{SOC}$ denotes the order of the cone $\mathcal{K}$.
    
    Recall the KKT conditions in \eqref{equ2: KKT conditions} and rewrite them into the following form:
    \begin{equation}
        \begin{split}
            \begin{bmatrix} 0 \\ 0 \\ s \end{bmatrix} = 
            \begin{bmatrix}
                Q & A^{T} & G^{T} \\
                -A & 0 & 0 \\
                -G & 0 & 0
            \end{bmatrix}
            \begin{bmatrix}
                x \\ y \\ z
            \end{bmatrix}
            + \begin{bmatrix} q \\ b \\ h \end{bmatrix} \\
            s^{T}z = 0, \quad (s, z) \in \mathcal{K} \times \mathcal{K}^{\ast}
        \end{split}
        \label{equ3: KKT conditions in matrix form}
    \end{equation}
    Then, solving \eqref{equ3: KKT conditions in matrix form} for $(x, y, z, s)$ is equivalent to solving the following monotone complementarity problem (MCP) over the cone $\mathcal{K}$ \cite{odonoghue2021operator, faybusovich1997euclidean}.
    \begin{equation}
        \begin{split}
            \underset{x, y, z, s}{\text{minimize}} \ & s^{T}z \\
            \text{subject to} \ & \begin{bmatrix} 0 \\ 0 \\ s \end{bmatrix} = 
            \begin{bmatrix}
                Q & A^{T} & G^{T} \\
                -A & 0 & 0 \\
                -G & 0 & 0
            \end{bmatrix}
            \begin{bmatrix}
                x \\ y \\ z
            \end{bmatrix}
            + \begin{bmatrix} q \\ b \\ h \end{bmatrix} \\
            & \ (s, z) \in \mathcal{K} \times \mathcal{K}^{\ast}
        \end{split}
        \tag{MCP}
        \label{equ3: MCP}
    \end{equation}
    \noindent
    Applying the barrier function to the \ref{equ3: MCP} yields the following optimization problem parameterized by the barrier parameter $\beta >0$.
    \begin{equation}
        \begin{split}
            \underset{x, y, z, s}{\text{minimize}} \quad & s^{T}z + \beta\Phi(s) + \beta\Phi(z) \\
            \text{subject to} \quad & \begin{bmatrix} 0 \\ 0 \\ s \end{bmatrix} = 
            \begin{bmatrix}
                Q & A^{T} & G^{T} \\
                -A & 0 & 0 \\
                -G & 0 & 0
            \end{bmatrix}
            \begin{bmatrix}
                x \\ y \\ z
            \end{bmatrix}
            + \begin{bmatrix} q \\ b \\ h \end{bmatrix} \\
            & (s, z) \in \text{int}\mathcal{K} \times \text{int}\mathcal{K}^{\ast}
        \end{split}
        \label{equ3: beta-parameterized problem}
    \end{equation}
    where the prefix $\text{int}(\cdot)$ represents the interior of the corresponding cone.
    The corresponding optimality conditions for the problem \eqref{equ3: beta-parameterized problem} are given by \cite{faybusovich1997euclidean}:

    \begin{equation}
        \begin{split}
            \begin{bmatrix} 0 \\ 0 \\ s \end{bmatrix} = 
            \begin{bmatrix}
                Q & A^{T} & G^{T} \\
                -A & 0 & 0 \\
                -G & 0 & 0
            \end{bmatrix}
            \begin{bmatrix}
                x \\ y \\ z
            \end{bmatrix}
            + \begin{bmatrix} q \\ b \\ h \end{bmatrix} \\
            s \circ z = \beta e,\quad (s, z) \in \text{int}\mathcal{K} \times \text{int}\mathcal{K}^{\ast}
        \end{split}
        \label{equ3: central path eq}
    \end{equation}
    \noindent
    where the operator $\circ$ represents Jordan product for Euclidean Jordan algebra and $e$ is the identity elements for Jordan product \cite{schmieta2003extension}. For every $\beta > 0$, there exists a unique solution $(x^{\ast}(\beta), y^{\ast}(\beta), z^{\ast}(\beta), s^{\ast}(\beta))$ to \eqref{equ3: central path eq}, called the \textit{central point}, and the set of all central points forms a smooth curve referred to as the \textit{central path} \cite{faybusovich1997euclidean}. As $\beta \rightarrow 0$, the central point $(x^{\ast}(\beta), y^{\ast}(\beta), z^{\ast}(\beta), s^{\ast}(\beta))$ converges to a solution of the MCP \cite{faybusovich1997euclidean}. Thus, the central path provides an effective guideline for the convergence of IPM iterates to the solution. In IPMs, the central path equations \eqref{equ3: central path eq} are solved by taking Newton steps while varying $\beta \rightarrow 0$ in each iteration.

\subsection{Homogeneous Embedding for Quadratic Cost Functions} \label{subsec : homogeneous embedding}
    
    Although IPMs are known to be efficient and robust for solving convex optimization problems, the standard implementation of IPMs has a drawback in detecting infeasibility of the problem \cite{meszaros2015practical}. To overcome this limitation, homogeneous embedding frameworks have been proposed. These reformulate the original problem into a lifted problem that exhibits the following desirable properties \cite{ye1994An, nesterov1999infeasible, andersen1999homogeneous, yoshise2007interior}: 
    \begin{enumerate}
        \item It has a bounded central path with a trivial initial point, regardless of the feasibility of the original problem.
        \item If the original problem is solvable, the solution of the lifted problem provides  certificates of optimality, and it can be mapped back to recover the solution of the original problem.
        \item If the original problem is infeasible, the solution of the lifted problem provides certificates of infeasibility.
    \end{enumerate}
    
    In this paper, we adopt the homogeneous embedding model for MCP (HMCP) over symmetric cone \cite{yoshise2007interior}, given by
    \begin{equation}
        \begin{split}
            &\underset{x, y, z, s, \kappa, \tau}{\text{minimize}} \quad s^{T}z + \kappa\tau \\
            &\text{subject to }\\
            & \ \begin{bmatrix} 0 \\ 0 \\ s \\ \kappa \end{bmatrix} = 
            \begin{bmatrix}
                Q & A^{T} & G^{T} & q \\
                -A & 0 & 0 & b \\
                -G & 0 & 0 & h \\
                -q^{T} & -b^{T} & -h^{T} & 0
            \end{bmatrix}
            \begin{bmatrix} x \\ y \\ z \\ \tau \end{bmatrix} + \begin{bmatrix} 0 \\ 0 \\ 0 \\ -\frac{1}{\tau} x^{T}Qx \end{bmatrix} \\
            & \ (s, z) \in \mathcal{K} \times \mathcal{K}^{\ast}, \quad (\kappa, \tau) \geq 0
        \end{split}
        \tag{HMCP}
        \label{equ3: HMCP}
    \end{equation}
    \noindent    
    The \ref{equ3: HMCP} remains feasible regardless of the feasibility of the original problem, and its solution $(x^{\ast}, y^{\ast}, z^{\ast}, s^{\ast}, \kappa^{\ast}, \tau^{\ast})$ provides the following certificates of optimality or infeasibility \cite{yoshise2007interior}:

    \begin{enumerate}
        \item \textit{Certificate of Optimality} : The \ref{equ3: MCP} has a solution if and only if $\tau^{\ast}>0$. In this case, $(x^{\ast}/\tau^{\ast}, s^{\ast}/\tau^{\ast})$ is the optimal solution of \ref{equ2: primal problem}, and $(y^{\ast}/\tau^{\ast}, z^{\ast}/\tau^{\ast})$ is the optimal solution of \ref{equ2: dual problem}.
        \item \textit{Certificate of Infeasibility} : The \ref{equ3: MCP} is infeasible if and only if $\kappa^{\ast} > 0$. In this case, the solution satisfies at least one of the following conditions:
        \begin{enumerate}
            \item $b^{T}y^{\ast} + h^{T}z^{\ast} < 0$ and the primal problem \ref{equ2: primal problem} is infeasible.
            \item $q^{T}x^{\ast} < 0$ and the dual problem \ref{equ2: dual problem} is infeasible.
        \end{enumerate}
    \end{enumerate}
    \noindent
    Leveraging these properties, we apply the IPM algorithm to the \ref{equ3: HMCP} instead of the original \ref{equ3: MCP}, ensuring convergence regardless of the feasibility of the original problem. The central path equations for \ref{equ3: HMCP} are given by:
    \begin{subequations}
        \begin{equation}
            \begin{bmatrix} 0 \\ 0 \\ s \\ \kappa \end{bmatrix} = 
            \begin{bmatrix}
                Q & A^{T} & G^{T} & q \\
                -A & 0 & 0 & b \\
                -G & 0 & 0 & h \\
                -q^{T} & -b^{T} & -h^{T} & 0
            \end{bmatrix}
            \begin{bmatrix} x \\ y \\ z \\ \tau \end{bmatrix} + \begin{bmatrix} 0 \\ 0 \\ 0 \\ -\frac{1}{\tau} x^{T}Qx \end{bmatrix} 
            \label{equ3: central path eq for HMCP - a}
        \end{equation}
        \begin{align}
            s \circ z = \beta e \label{equ3: central path eq for HMCP - b} \\
            \kappa \tau = \beta \label{equ3: central path eq for HMCP - c} \\
            (s, z) \in \text{int}\mathcal{K} \times \text{int}\mathcal{K}^{\ast}, \quad (\kappa, \tau) > 0 \label{equ3: central path eq for HMCP - d}
        \end{align}
        \label{equ3: central path eq for HMCP}
    \end{subequations}
    Note that the quadratic objective term is explicitly accounted for in the formulation.
    As a result, unlike the conventional homogeneous self-dual embedding \cite{ye1994An, nesterov1999infeasible}, the presented \ref{equ3: HMCP} does not require problem reformulation that introduces an additional slack variable $r$ and a SOC constraint to replace the quadratic objective term, as in \eqref{equ3: linear conic program} \cite{lobo2998applications}.
    This prevents unnecessary degradation in solver performance caused by increased problem size and reduced sparsity.
    Moreover, since the \ref{equ3: HMCP} becomes identical to the homogeneous self-dual embedding when $Q=0$, it can be regarded as a more general homogeneous embedding formulation.

    \begin{equation}
        \begin{split}
            \underset{x, s, r}{\text{minimize}} \quad & q^{T}x + r \\
            \text{subject to} \quad & Ax = b, \quad Gx + s = h \\
            & \left\Vert \begin{matrix} r-\frac{1}{2} \\ Q^{1/2}x  \end{matrix} \right\Vert_{2} \leq r + \frac{1}{2} \\
            & s \in \mathcal{K}, \quad r \geq 0
        \end{split}
        \label{equ3: linear conic program}
    \end{equation}
    \noindent

\subsection{Nesterov-Todd Scaling}

    Nesterov and Todd extended PDIPMs to a more general setting than linear programming by developing the concepts of self-scaled barriers and self-scaled cones \cite{nesterov1997self-scaled, nesterov1998primal-dual}. The cone $\mathcal{K}$ in this paper, given by \eqref{equ2 : cone K}, also belongs to the class of self-scaled cones. In fact, it has been shown that the family of self-scaled cones is identical to the set of symmetric cones \cite{guler1996barrier}. Nesterov and Todd showed that there exists a unique scaling point $w \in \text{int}\mathcal{K}$ corresponding to any points $s \in \text{int}\mathcal{K}$ and $z \in \text{int}\mathcal{K}^{\ast}$ such that
    \begin{equation}
        H(w)s = z
        \label{equ3: self-scaled properties}
    \end{equation} 
    where $H$ represents the Hessian of the barrier function $\Phi$ \cite{nesterov1997self-scaled}. This indicates that there exists a class of scaling linear transformations for such cones, which come from the Hessians of the barrier evaluated at the points within the cone itself, i.e., self-scaled. Using the properties of the Hessian of the log-det barrier \cite{vandenberghe2010cvxopt}, the scaling point $w$ is given by:
    \begin{equation}
        w  = H( s^{1/2})^{-1} \left( H(s^{1/2})^{-1}z \right)^{-1/2}
        \label{equ3: scaling point}
    \end{equation}
    Then, the Nesterov-Todd scaling (NT-scaling) matrix $W$ is obtained by factoring $H(w)$ as:
    \begin{equation}
        H(w)^{-1} = W^{T}W
        \label{equ3: factoring hessian of barrier}
    \end{equation}
    The symmetric scaling matrix can be chosen as $W = H(w)^{-1/2}$, since $H(w)$ is positive definite (and thus $W$ is) \cite{nesterov1997self-scaled}. Notably, this scaling leaves the interior of the cone and the central path invariant \cite{vandenberghe2010cvxopt}, i.e., 
    \begin{equation*}
        s \in \text{int}\mathcal{K} \Leftrightarrow W^{-T}s \in \text{int}\mathcal{K}, \quad z \in \mathcal{K} \Leftrightarrow Wz \in \text{int}\mathcal{K}
    \end{equation*}
    \begin{equation*}
        s \circ z  = \beta e \ \Leftrightarrow \ (W^{-T}s) \circ (Wz) = \beta e
    \end{equation*}
    \noindent
    Thus, the scaled central path equations can be constructed from \eqref{equ3: central path eq for HMCP}, while replacing \eqref{equ3: central path eq for HMCP - b} with \eqref{equ3: scaled complementary slackness}.
    \begin{equation}
        (W^{-T}s)\circ(Wz) = \beta e
        \label{equ3: scaled complementary slackness}
    \end{equation}
    \noindent
    From \eqref{equ3: self-scaled properties} and \eqref{equ3: factoring hessian of barrier}, it can also be shown that the scaling maps the points $s$, $z$ to the same scaled point $\lambda$. 
    \begin{equation}
        \lambda = W^{-T}s = Wz
        \label{equ3: scaled point}
    \end{equation}
    \noindent
    Lastly, we provide the analytic expressions of the NT-scaling components for the cones of interest in this paper \cite{vandenberghe2010cvxopt}: 
    \begin{enumerate}
        \item \textit{NNC} :
        For $s, \ z \in \text{int}\mathbb{R}_{+}$, 
        \begin{itemize}
            \begin{subequations}
            \item Scaling Point ($w$) :  
                \begin{equation}
                    w = s^{1/2} \circ z^{-1/2} = \sqrt{\frac{s}{z}}
                \end{equation}
            \item Scaling Matrix ($W$) :
                \begin{equation}
                    W = \text{diag}(w) = \sqrt{\frac{s}{z}}
                \end{equation}
            \item Scaled Point ($\lambda$) :
                \begin{equation}
                    \lambda = s^{1/2} \circ z^{1/2} = \sqrt{s \cdot z}
                \end{equation}
            \label{equ3: NNC scaling}
            \end{subequations}
        \end{itemize}
        \item \textit{SOC} : 
        For $s, \ z \in \text{int}\mathbb{Q}^{d}$ with $d >  1$, let $\bar{s} = (\bar{s}_{0}, \bar{s}_{1})$, $\bar{z} = (\bar{z}_{0}, \bar{z}_{1}) \in \mathbb{R} \times \mathbb{R}^{d-1}$ be the normalized points
            \begin{equation*}
                \bar{s} = \frac{1}{\left( s_{0}^{2} - s_{1}^{T}s_{1} \right)^{1/2}} s, \quad \bar{z} = \frac{1}{\left( z_{0}^{2} - z_{1}^{T}z_{1} \right)^{1/2}} z
            \end{equation*}
            and $\gamma$, $\bar{w} = (\bar{w}_{0}, \bar{w}_{1}) \in \mathbb{R} \times \mathbb{R}^{d-1}$, $\eta$ be
            \begin{equation}
                \gamma = \left( \frac{1 + \bar{s}^{T}\bar{z}}{2} \right)^{1/2}
            \end{equation}
            \begin{equation}
                \bar{w} = \frac{1}{2\gamma} \left( \bar{s} + \begin{bmatrix} \bar{z}_{0} \\ -\bar{z}_{1} \end{bmatrix} \right)
            \end{equation}
            \begin{equation}
                \eta = \left( \frac{s_{0}^{2} - s_{1}^{T}s_{1}}{z_{0}^{2} - z_{1}^{T}z_{1}} \right)^{1/4} 
            \end{equation}
        \begin{itemize}
            \begin{subequations}
            \item Scaling Point ($w$) :
                \begin{equation}
                            w = \eta \bar{w}
                \end{equation}
            \item Scaling Matrix ($W$) :
                \begin{equation}
                    W = \eta
                    \begin{bmatrix}
                        \bar{w}_{0} & \bar{w}_{1}^{T} \\
                        \bar{w}_{1} & I_{d-1} + (1 + \bar{w}_{0})^{-1}\bar{w}_{1} \bar{w}_{1}^{T}
                    \end{bmatrix}
                    \label{equ3: Nesterov-Todd scaling for SOC}
                \end{equation}
            \item Scaled Point ($\lambda$) :
                \begin{equation}
                    \lambda = \eta
                    \begin{bmatrix}
                        \bar{w}_{0}z_{0} + \bar{w}_{1}^{T}z_{1} \\
                        z_{1} + \displaystyle \left( z_{0} + \frac{\bar{w}_{1}^{T}z_{1}}{1 + \bar{w}_{0}} \right) \bar{w}_{1}
                    \end{bmatrix}
                \end{equation}
                \label{equ3: SOC scaling}
            \end{subequations}
        \end{itemize}
    \end{enumerate}
    Note that $1+\bar{w}_{0} \neq 0$ since $\bar{w}_{0} > 0$ from its definition.
    
    For the cone $\mathcal{K}$ in \eqref{equ2 : cone K}, the NT-scaling matrix is constructed as a block diagonal matrix composed of scaling matrices for each NNC and SOC.
    \begin{equation}
        W = \texttt{blkdiag}(W_{1}, \cdots, W_{l+n_{SOC}} )
    \end{equation}
    where the first $l$ blocks correspond to the scaling matrices for NNCs and the remaining $n_{SOC}$ blocks correspond to those for SOCs.

\subsection{Search Directions} \label{subsec: search directions}

    The solver developed in this paper utilizes PDIPMs, which are known to be among the most efficient frameworks within IPMs \cite{nocedal1999numerical}. PDIPMs compute a Newton-like search direction by linearizing the (scaled) central path equations and loosely follow the central path toward the solution. 
    Since the intermediate iterates of PDIPMs do not lie exactly on the central path, they are not generally primal-dual feasible to \eqref{equ2: primal problem} and \eqref{equ2: dual problem}, as well as equations \eqref{equ3: central path eq for HMCP - b}, \eqref{equ3: central path eq for HMCP - c} or \eqref{equ3: scaled complementary slackness} are not precisely satisfied for the barrier parameter $\beta$ \cite[Chapter 10]{wright1997primal}. Therefore, duality gap cannot be evaluated as $s^{T}z+\kappa\tau$ at intermediate iterates. Instead, PDIPMs introduce the \textit{duality measure} $\mu$ as an indicator of the approximated duality gap at each point in the search space instead of $\beta$ \cite{wright1997primal, nocedal1999numerical}.
    \begin{equation}
        \mu = \frac{s^{T}z + \kappa\tau}{l+n_{SOC}+1} >0
        \label{equ3: mu}
    \end{equation}
    
    For the iterates of PDIPMs to converge stably to the solution, appropriate search directions should be determined by balancing the trade-off between the following two objectives: 1) reducing the duality gap to achieve optimality, and 2) tracing the central path to maintain feasibility. To achieve this, PDIPMs modify the Newton step by biasing it toward a point satisfying
    \begin{equation}
        (W^{-T}s) \circ (Wz) = \sigma \mu e, \quad \kappa \tau = \sigma \mu e
        \label{equ3: centered scaled complementarity}
    \end{equation}
    where $\sigma \in [0, 1]$ is referred to as a \textit{centering parameter}, which determines the degree of bias in the search direction. If $\sigma = 0$, the search direction is called the \textit{affine search direction}, which corresponds to a pure Newton step aimed at making progress toward optimality. If $\sigma = 1$, it is called the \textit{centering direction} and is directed toward the central point corresponding to the current duality measure $\mu$. Although the centering direction usually makes little progress in reducing the duality gap, it enables greater progress in the next step by bringing the iterates close to the central path \cite{wright1997primal}.

    Mehrotra proposed a practically effective predictor-corrector method to compute search directions based on heuristics for selecting the centering parameter \cite{mehrotra1992on}. By linearizing the scaled central path equations and introducing the centering parameter, we obtain the following linear equations for computing search directions.
    \begin{subequations}
        \begin{align}
            \begin{bmatrix} 0 \\ 0 \\ \Delta s \\ \Delta \kappa \end{bmatrix} - 
            \Lambda \begin{bmatrix} \Delta x \\ \Delta y \\ \Delta z \\ \Delta \tau \end{bmatrix} = -(1-\sigma) 
            \begin{bmatrix} r_{x} \\ r_{y} \\ r_{z} \\ r_{\tau} \end{bmatrix} \\
            \lambda \circ \left( W^{-T}\Delta s + W \Delta z \right) = -r_{s} \label{equ3: search direction eq. - b} \\
            \kappa \Delta \tau + \tau \Delta \kappa = -r_{\kappa} \label{equ3: search direction eq. - c}
        \end{align}
        \label{equ3: search direction eq.}
    \end{subequations}
    \noindent
    where
    \begin{equation}
        \Lambda := \begin{bmatrix}
                Q & A^{T} & G^{T} & q \\
                -A & 0 & 0 & b \\
                -G & 0 & 0 & h \\
                -q^{T} - \frac{2}{\tau}x^{T}Q & -b^{T} & -h^{T} & \frac{1}{\tau^{2}}x^{T}Qx
            \end{bmatrix}
        \label{equ3: coefficient matrix Lambda}
    \end{equation}
    and $(r_{x}, r_{y}, r_{z}, r_{\tau}, r_{s}, r_{\kappa})$ are the residuals of the central path equations for the current iterates $\chi = (x, y, z, s, \kappa, \tau)$.
    \begin{equation}
        \begin{split}
            r_{x} &= -Qx -A^{T}y - G^{T}z - q\tau \\
            r_{y} &= Ax - b\tau \\
            r_{z} &= Gx + s -h\tau \\
            r_{\tau} &= q^{T}x + b^{T}y + h^{T}z + \frac{1}{\tau}x^{T}Qx + \kappa
        \end{split}
        \label{equ3: residuals}
    \end{equation}
    \begin{equation}
        \begin{split}
            &r_{s} = \lambda \circ \lambda - \sigma \mu e \qquad \qquad \qquad \qquad \quad \\
            &r_{\kappa} = \kappa \tau - \sigma \mu 
        \end{split}
        \label{equ3: residuals 2}
    \end{equation}
    \noindent
    Mehrotra's predictor-corrector methods utilize the affine search direction as a predictor to determine the level of centering and to estimate the complementarity error arising from linearization. The affine search direction is computed by solving \eqref{equ3: search direction eq.} with $\sigma = 0$ and is denoted as $\Delta \chi_{a} = \Delta(x_{a}, y_{a}, z_{a}, s_{a}, \kappa_{a}, \tau_{a})$. Based on Mehrotra's heuristics, the centering parameter is computed as
    \begin{equation}
        \sigma = \left(\frac{\mu_{a}}{\mu}\right)^{3} = (1- \alpha_{a})^{3}
        \label{equ3: centering parameter}
    \end{equation}
    where $\mu_{a}$ is the predicted duality measure for the affine search direction
    \begin{equation}
        \mu_{a} = \frac{(s+\alpha_{a}\Delta s_{a})^{T}(z + \alpha_{a}\Delta z_{a}) + (\kappa + \alpha_{a} \Delta \kappa_{a})(\tau + \alpha_{a} \Delta \tau_{a})}{l+n_{SOC}+1}
        \label{equ3: mu_a}
    \end{equation}
    and $\alpha_{a}$ is the step size obtained from the line search in \eqref{equ3: linesearch} along the affine search direction (i.e., by replacing $\Delta(s, z, \kappa, \tau)$ in \eqref{equ3: linesearch} with $\Delta(s_{a}, z_{a}, \kappa_{a}, \tau_{a})$).
    \begin{equation}
        \begin{split}
            \alpha_{\text{step}} = \sup \{\bar{\alpha} \in (0, 1) \left. \right\vert \ & (s, z) + \bar{\alpha} (\Delta s, \Delta z) \in  \text{int}\mathcal{K} \times \text{int}\mathcal{K}^{\ast}, \\ &(\kappa, \tau) + \bar{\alpha} (\Delta \kappa, \Delta \tau) > 0 \}
            \label{equ3: linesearch}
        \end{split}
    \end{equation}
    Furthermore, the nonlinearity errors in complementarity associated with the affine search direction are computed as follows:
    \begin{subequations}
        \begin{align}
            \begin{split}
                e_{sz} :&= (W^{-T}(s + \Delta s_{a}))\circ(W(z+\Delta z_{a})) - 0 \\
                & = (W^{-T}\Delta s_{a}) \circ (W \Delta z_{a})
            \end{split}
            \\
            \begin{split}
                e_{\kappa\tau} :&= (\kappa + \Delta\kappa_{a})(\tau + \Delta \tau_{a}) - 0 \\
                &= \Delta \kappa_{a} \Delta \tau_{a}
            \end{split}
        \end{align}
        \label{equ3: nonlinearity error}
    \end{subequations}
    
    In the corrector step, the nonlinearity errors are accounted for by replacing $r_{s}$ and $r_{\kappa}$ in \eqref{equ3: residuals 2} with those given in \eqref{equ3: rhs of nonlinearity corrected complementarity}.
    \begin{subequations}
        \begin{align}
            r_{s} =& \lambda \circ \lambda - \sigma \mu e + (W^{-T}\Delta s_{a})\circ(W \Delta z_{a}) \\
            r_{\kappa} = & \kappa \tau - \sigma \mu +\Delta \kappa_{a} \Delta \tau_{a}
        \end{align}
        \label{equ3: rhs of nonlinearity corrected complementarity}
    \end{subequations}
    Then, the search direction is computed by solving \eqref{equ3: search direction eq.} with the centering parameter from \eqref{equ3: centering parameter} and the residuals in \eqref{equ3: rhs of nonlinearity corrected complementarity}.
    We refer to this search direction as \textit{combined search direction} and denote it as $\Delta\chi_{c} = \Delta(x_{c}, y_{c}, z_{c}, s_{c}, \kappa_{c}, \tau_{c})$. Finally, the iterates are updated by taking a step along the combined search direction 
    \begin{equation}
        \chi \leftarrow \chi + \alpha \Delta \chi_{c}
        \label{equ3: iterates update}
    \end{equation}
    where the step size $\alpha$ is obtained from the line search in \eqref{equ3: linesearch} along the combined search direction $\Delta (s_{c}, z_{c}, \kappa_{c}, \tau_{c})$.
    
    For both the affine and combined search directions, we need to solve the following form of the linear system with the corresponding right-hand sides $(d_{x}, d_{y}, d_{z}, d_{\tau}, d_{s}, d_{\kappa})$.
    \begin{subequations}
        \begin{align}
             \begin{bmatrix} 0 \\ 0 \\ \Delta s \\ \Delta \kappa \end{bmatrix} - 
            \Lambda \begin{bmatrix} \Delta x \\ \Delta y \\ \Delta z \\ \Delta \tau \end{bmatrix} = 
            \begin{bmatrix} d_{x} \\ d_{y} \\ d_{z} \\ d_{\tau} \end{bmatrix} \label{equ3: linear system 1} \\
            \lambda \circ \left( W^{-T}\Delta s + W \Delta z \right) = d_{s} \label{equ3: linear system 2} \\
            \kappa \Delta \tau + \tau \Delta \kappa = d_{\kappa}  \label{equ3: linear system 3}
        \end{align}
        \label{equ3: linear system}
    \end{subequations}
    Solving this linear system is the most computationally expensive step in PDIPMs, thereby determining the performance of the solver. Therefore, we adopt the approach of \cite{vandenberghe2010cvxopt} to efficiently compute the search directions. To simplify notation, we define the following vectors
    \begin{subequations}
        \begin{align}
            \Delta \xi &:= \begin{bmatrix} \Delta x^{T} & \Delta y^{T} & \Delta z^{T} \end{bmatrix}^{T} \\
            b_{1} &:= \begin{bmatrix} -q^{T} & b^{T} & h^{T} \end{bmatrix}^{T} \\
            b_{2} &:= \begin{bmatrix} -d_{x}^{T} & d_{y}^{T} & (d_{z} - W(\lambda \backslash d_{s}))^{T} \end{bmatrix}^{T}
        \end{align}
    \end{subequations}
    where the operator $\backslash$ is the inverse operator to the Jordan product $\circ$ such that if $u \circ v = w$, then $u \backslash w = v$. (Refer to \cite{domahidi2013ECOS} for the analytical formula of the inverse operator for the NNCs and SOCs.) We also define the coefficient matrix $K$, referred to as the \textit{KKT matrix}, as follows:
    \begin{equation}
        K := \begin{bmatrix}
            Q & A^{T} & G^{T} \\
            A & 0 & 0 \\
            G & 0 & -W^{T}W
        \end{bmatrix}
        \label{equ3: kkt matrix}
    \end{equation}
    By eliminating $\Delta s$ and $\Delta \kappa$ from \eqref{equ3: linear system 1} using \eqref{equ3: linear system 2} and \eqref{equ3: linear system 3}, the linear system \eqref{equ3: linear system} can be reformulated as follows:
    \begin{subequations}
        \begin{equation}
             K \Delta \xi - b_{1} \Delta \tau  = b_{2}
            \label{equ3: eliminated linear system 1}
        \end{equation}
        \begin{equation}
            \Delta \tau = \frac{-d_{\tau} + d_{\kappa}/\tau + \left( q^{T} + \frac{2}{\tau}x^{T}Q \right)\Delta x + b^{T}\Delta y + h^{T} \Delta z}{\frac{1}{\tau^{2}}x^{T}Qx + \kappa/\tau}
            \label{equ3: eliminated linear system 2}
        \end{equation}
    \end{subequations}    
    The linear system in \eqref{equ3: eliminated linear system 1} can be solved using a pair of linear systems, which share a common coefficient matrix $K$.
    \begin{subequations}
        \begin{align}
            K\Delta \xi_{1}  = b_{1} \label{equ3: KKT system 1} \\
            K\Delta \xi_{2} = b_{2} \label{equ3: KKT system 2}
        \end{align}
        \label{equ3: linear system pair}
    \end{subequations}
    By solving \eqref{equ3: linear system pair} for $\Delta \xi_{1}$ and $\Delta \xi_{2}$, the search directions are recovered as a linear combination of $\Delta \xi_{1}$ and $\Delta \xi_{2}$:
    \begin{equation}
        \Delta \xi  = \Delta \xi_{2} + \Delta \tau \Delta \xi_{1}
    \end{equation}
    Subsequently, $\Delta s$ and $\Delta \kappa$ can be recovered from \eqref{equ3: linear system 2} and \eqref{equ3: linear system 3}, as follows:

    \begin{subequations}
        \begin{align}
            \Delta s &= W ( \lambda \backslash d_{s} - W \Delta z) \\
            \Delta \kappa &= (d_{\kappa} - \kappa \Delta \tau) / \tau
        \end{align}
    \end{subequations}

    It is important to note that solving \eqref{equ3: KKT system 1} is the same for both the affine and combined search direction. Therefore, it only needs to be solved once in each iteration. In contrast, the equation \eqref{equ3: KKT system 2} has to be solved separately for each search direction with the corresponding $b_{2}$. Consequently, the linear systems need to be solved three times in total to compute the search direction. However, since all three linear systems share the same coefficient matrix $K$, the matrix factorization—which requires the most computational cost—needs to be performed only once. For the different right-hand sides, only forward and backward substitutions are required, and their computational costs are negligible compared to that of the matrix factorization.

    To further improve computational efficiency in matrix factorization, we adopt the structure exploitation strategy for SOC scaling proposed in \cite{domahidi2013ECOS} to obtain a sparse KKT matrix. For a SOC $\mathbb{Q}^{d}$, the scaling matrix $W_{SOC}$ in \eqref{equ3: Nesterov-Todd scaling for SOC} propagates to the (3, 3) block in the KKT matrix as $-W^{T}_{SOC}W_{SOC}$, which is a dense matrix with $d^{2}$ nonzero entries. This results in a significant number of nonzeros in both the KKT matrix and its factorization, even for problems with sparse data. By exploiting the diagonal-plus-rank-one structure of $W_{SOC}$, we can expand $W^{T}_{SOC}W_{SOC}$ to a sparse form as
    \begin{equation}
        \tilde{W}_{SOC}^{T}\tilde{W}_{SOC} = \eta^{2} 
        \begin{bmatrix}
            \bar{D} & v & u \\
            v^{T} & 1 & 0 \\
            u^{T} & 0 & -1
        \end{bmatrix}
    \end{equation}
    where
    \begin{align*}
        \bar{D} := \begin{bmatrix} a & \\ & I_{d-1} \end{bmatrix}, \ u := \begin{bmatrix} u_{0} \\ u_{1}\bar{w}_{1} \end{bmatrix}, \ v := \begin{bmatrix} 0 \\ v_{1}\bar{w}_{1} \end{bmatrix} \\
        u_{0} := \sqrt{\bar{w}_{0}^{2} + \Vert \bar{w}_{1} \Vert_{2}^{2}-a}, \ u_{1} := \frac{c}{u_{0}}, \ v_{1} := \sqrt{\frac{c^{2}}{u_{0}^{2}}-\bar{d}} \\
        a := \frac{1}{2}\left( \bar{w}_{0}^{2} + \Vert \bar{w}_{1} \Vert^{2} - \frac{c^{2}\Vert \bar{w} \Vert_{2}^{2}}{1 + d \Vert \bar{w} \Vert_{2}^{2}} \right)\\
        c := 1 + \bar{w}_{0} + \frac{\Vert \bar{w}_{1} \Vert_{2}^{2}}{1 + \bar{w}_{0}}, \ \bar{d} := 1 + \frac{2}{1 + \bar{w}_{0}} + \frac{\Vert \bar{w}_{1} \Vert_{2}^{2}}{(1 + \bar{w}_{0})^{2}}
    \end{align*}

\subsection{Solving the Linear System} \label{subsec: solving the linear system}

    Computing the search directions requires solving the following form of the linear system.
    \begin{equation}
        K \Delta \xi_{i} = b_{i}
        \label{equ3: KKT system form}
    \end{equation}
    where $i = 1, \ 2$.
    Since the KKT matrix $K$ is symmetric indefinite, we apply \textit{static regularization} \cite{vanderbei1995symmetric, mattingley2012cvxgen} to render the KKT matrix \textit{quasi-definite},
    i.e., symmetric with the block located in the (1, 1) position being positive definite, the block in the (2, 2) position being negative definite, and the off-diagonal blocks being full-rank \cite{vanderbei1995symmetric}.
    \begin{equation}
        \tilde{K} = \left[ \begin{array}{c|cc} Q & A^{T} & G^{T} \\
            \hline
            A & 0 & 0 \\
            G & 0 & -W^{T}W  \end{array} \right] + 
            \left[ \begin{array}{c|c} \delta_{s} I & 0 \\ 
            \hline
            0 & -\delta_{s}I
            \end{array} \right]
    \end{equation}
    where $\delta_{s} > 0$ is a static regularization parameter and the value of $10^{-7}$ is used as a default in this implementation. For such matrix $\tilde{K}$, a stable $LDL^{T}$ factorization
    \begin{equation}
        P\tilde{K}P^{T} = LDL^{T}
    \end{equation}
    exists for any permutation $P$, where $L$ denotes a lower triangular matrix and $D$ a diagonal matrix in the resulting factorization. This property eliminates the need for data-dependent permutations, which are typically required in the $LDL^{T}$ factorization of general symmetric indefinite matrices \cite{bunch1971direct}. As a result, the process of solving the linear system—which requires the most computational effort—can be formulated more efficiently.

    A commonly used strategy for selecting the permutation $P$ is to minimize \textit{fill-in} in the factor $L$. Fill-in refers to additional nonzero entries that arise in $L$ during the factorization but are not present in the original matrix. These entries increase memory usage and the number of floating-point operations, leading to a degradation in solver performance. However, finding an optimal permutation is known to be NP-complete \cite{yannakakis1981computing}. Therefore, a heuristic method known as approximate minimum degree (AMD) ordering \cite{amestoy1996approximate} is employed to compute a permutation that effectively reduces fill-in. Moreover, since AMD ordering depends only on the sparsity pattern (the location of the nonzero elements) of the matrix to be factored, it needs to be performed only once at the initial stage of the algorithm.

    Although the $LDL^{T}$ factorization is theoretically stable for $\tilde{K}$, it is still necessary to ensure numerical stability in practice due to round-off errors introduced by the finite precision of floating-point arithmetic. To this end, \textit{dynamic regularization} is applied during the $LDL^{T}$ factorization process. Specifically, it modifies the entries of the diagonal factor $D$ whenever they become close to zero or change sign:
    \begin{equation*}
        \text{If } \vert D_{i}\vert \leq \epsilon_{d} \quad \Rightarrow \quad D_{i} \leftarrow S_{i}\delta_{d}
    \end{equation*}
    where $D_{i}$ is the $i$th diagonal entry of $D$, and $S_{i}$ is the sign of $D_{i}$ which can be predetermined from the sign of the corresponding diagonal entry in the KKT matrix. For $\epsilon_{d} > 0$ and $\delta_{d} > 0$, we set $\epsilon_{d} = 10^{-13}$ and $\delta_{d} = 10^{-7}$ in this implementation.

    After obtaining the factors $L$ and $D$ through factorization, the linear system in \eqref{equ3: KKT system form} can be easily solved via backward substitution, diagonal scaling, and forward substitution \cite[\S3]{davis2006direct}.
    However, the resulting solution corresponds to the regularized system $\tilde{K}\Delta \xi_{i} = b_{i}$, rather than the original system $K \Delta \xi_{i} = b_{i}$. To recover the solution to the original system, we apply \textit{iterative refinement} \cite[\S4.11]{duff2017direct} to progressively reduce the residual error $\Vert K \Delta \xi_{i} - b_{i} \Vert_{\infty}$ in $\Delta \xi_{i}$. Specifically, it iteratively solves the linear system 
    \begin{equation}
        \tilde{K} \delta\xi_{i} = b_{i} -  K \Delta \xi_{i}
        \label{equ3: linear system IR}
    \end{equation}
    for $\delta \xi_{i}$, and updates $\Delta \xi_{i}$ as 
    \begin{equation}
        \Delta \xi_{i} \leftarrow \Delta \xi_{i} + \delta \xi_{i}
    \end{equation}
    until $\Vert K \Delta \xi_{i} - b_{i} \Vert_{\infty} < \epsilon_{ir}$, where the default tolerance is set to $\epsilon_{ir}=10^{-13}$.
    It is worth noting that no additional factorization is required to solve the linear system in \eqref{equ3: linear system IR}, as it shares the same coefficient matrix $\tilde{K}$ with the system in \eqref{equ3: KKT system form}. Therefore, only forward, backward substitution and diagonal scaling are needed during the iterative refinement process.
    
\subsection{Initialization} \label{subsec : initialization}

    For choosing an initial iterate of the proposed solver, we adopt the initialization strategies of \cite{vandenberghe2010cvxopt}. If $Q \neq 0$, the initial iterate $(x^{(0)}, y^{(0)}, z^{(0)}, s^{(0)})$ is chosen as the solution of the following constrained least-squares problem:
    \begin{equation}
        \begin{split}
            \underset{\hat{x}, \hat{s}}{\text{minimize}} & \quad \frac{1}{2}\hat{x}^{T}Q\hat{x} + q^{T}\hat{x} + \frac{1}{2}\Vert \hat{s} \Vert_{2}^{2} \\
            \text{subject to} & \quad G\hat{x} + \hat{s} = h, \quad A\hat{x} = b
        \end{split}
    \label{equ3: init problem (Q != 0)}
    \end{equation}
    where $\hat{(\cdot)}$ denotes the auxiliary variables associated with the initialization problem.
    Applying the optimality condition, the solution of the problem can be obtained by solving the following linear system:
    \begin{equation}
        \begin{bmatrix}
            Q & A^{T} & G^{T} \\
            A & 0 & 0 \\
            G & 0 & -I
        \end{bmatrix} \begin{bmatrix} \hat{x} \\ \hat{y} \\ \hat{z} \end{bmatrix} = \begin{bmatrix} -q \\ b \\ h \end{bmatrix}
    \label{equ3: linear system for Q != 0}
    \end{equation}
    We set $x^{(0)} = \hat{x}$, $y^{(0)} = \hat{y}$ and
    \begin{equation}
        s^{(0)} = 
        \begin{cases}
            -\hat{z} & \qquad \alpha_{p} < 0 \\
            -\hat{z} + (1 + \alpha_{p}) e & \qquad \text{otherwise}
        \end{cases}
    \end{equation}
    where $\alpha_{p} = \inf \{ \alpha : -\hat{z} + \alpha e \in \text{int}\mathcal{K} \}$ for $s^{(0)}$ to be located in the interior of the cone $\mathcal{K}$. Similarly for $z^{(0)}$, we set
    \begin{equation}
        z^{(0)} = 
        \begin{cases}
            \hat{z} & \qquad \alpha_{d} < 0 \\
            \hat{z} + (1 + \alpha_{d}) e & \qquad \text{otherwise}
        \end{cases}
    \end{equation}
    where $\alpha_{d} = \inf \{ \alpha : \hat{z} + \alpha e \in \text{int}\mathcal{K} \}$.
    
    If $Q = 0$, the initial primal variables $(x^{(0)}, s^{(0)})$ are obtained from the problem:
    \begin{equation}
        \begin{split}
            \underset{\hat{x}, \hat{s}}{\text{minimize}} & \quad \frac{1}{2} \Vert \hat{s} \Vert_{2}^{2} \\
            \text{subject to} & \quad A\hat{x} = b, \quad G\hat{x} + \hat{s} = h
        \end{split}
    \label{equ3: init probelm form primal (Q=0)}
    \end{equation}
    which provides feasible $x^{(0)}$ with respect to the equality constraints, and minimum norm of the residual for the inequality constraints. The solution of the problem \eqref{equ3: init probelm form primal (Q=0)} is obtained by solving the following linear system:
    \begin{equation}
        \begin{bmatrix}
            0 & A^{T} & G^{T} \\
            A & 0 & 0 \\
            G & 0 & -I
        \end{bmatrix} \begin{bmatrix} \hat{x} \\ \hat{y} \\ -\hat{s} \end{bmatrix} = 
        \begin{bmatrix}
            0 \\ b \\ h
        \end{bmatrix}
    \label{equ3: linear system for Q = 0, primal}
    \end{equation}
    We set $x^{(0)} = \hat{x}$, and
    \begin{equation}
        s^{(0)} = \begin{cases} \hat{s} &\qquad \alpha_{p} < 0 \\ \hat{s} + (1 + \alpha_{p})e &\qquad \text{otherwise}  \end{cases}
    \end{equation}
    where $\alpha_{p} = \inf \{ \alpha : \hat{s} + \alpha e \in \text{int}\mathcal{K} \}$. For initial dual variable $(y^{(0)}, z^{(0)})$, we solve the constrained least-squares problem to minimize the norm of the residual for the inequality constraints in the dual problem.
    \begin{equation}
        \begin{split}
            \underset{\hat{y}, \hat{z}}{\text{minimize}} & \quad \frac{1}{2} \Vert \hat{z} \Vert_{2}^{2} \\
            \text{subject to} & \quad A^{T}\hat{y} + G^{T}\hat{z} + q = 0
        \end{split}
    \label{equ3: init probelm form dual (Q=0)}
    \end{equation}
     The solution of the problem \eqref{equ3: init probelm form dual (Q=0)} can be obtained by solving the following linear system:
    \begin{equation}
        \begin{bmatrix}
            0 & A^{T} & G^{T} \\
            A & 0 & 0 \\
            G & 0 & -I
        \end{bmatrix} \begin{bmatrix} \hat{x} \\ \hat{y} \\ \hat{z} \end{bmatrix} =
        \begin{bmatrix}
            -q \\ 0 \\ 0
        \end{bmatrix}
    \label{equ3: linear system for Q = 0, dual}
    \end{equation}
    We set $y^{(0)} = \hat{y}$, and
    \begin{equation}
        z^{(0)} = \begin{cases} \hat{z} & \qquad \alpha_{d} < 0  \\ \hat{z} + (1 + \alpha_{d})e & \qquad \text{otherwise}  \end{cases}
    \end{equation}
    where $\alpha_{d} = \inf \{ \alpha : \hat{z} + \alpha e \in \text{int}\mathcal{K} \}$. We set $\kappa^{(0)} = 1$, $\tau^{(0)} = 1$ for both case of $Q = 0$ and $Q \neq 0$.

    This initialization process provides feasible points with respect to the equality constraints, and minimizes the norm of the residual for the inequality constraints, which helps prevent $\alpha_{p}$ and $\alpha_{d}$ from being computed excessively large. Moreover, it is important to note that the linear systems for the initialization share the same structure as those used for computing the search directions. This enables the construction of an efficient algorithm by allowing the use of the same sparsity pattern and permutation for the matrix factorization.

\subsection{Termination Criteria}

    The termination criteria of the solver are based on the certificates of optimality and infeasibility provided by the homogeneous embedding in \S \ref{subsec : homogeneous embedding}. The solver declares the convergence to the optimal solution if the following conditions are satisfied:
    \begin{subequations}
        \begin{equation}
            \frac{\Vert A(x/\tau) - b \Vert_{2}}{\max\{1, \ \Vert x \Vert_{2} + \Vert b \Vert_{2} \}} < \epsilon_{\text{feas}}
            \label{equ3: pres_eq}
        \end{equation}
        \begin{equation}
            \frac{\Vert G(x/\tau) + (s/\tau) - h\Vert_{2}}{\max\{1, \ \Vert x \Vert_{2} + \Vert s \Vert_{2} + \Vert h \Vert_{2} \}} < \epsilon_{\text{feas}}
            \label{equ3: pres_ineq}
        \end{equation}
        \begin{equation}
            \frac{\Vert Q(x/\tau) + A^{T}(y/\tau) + G^{T}(z/\tau) + q \Vert_{2}}{\max\{1, \ \Vert x \Vert_{2} + \Vert y \Vert_{2} + \Vert z \Vert_{2} + \Vert q \Vert_{2} \}} < \epsilon_{\text{feas}}
            \label{equ3: dres}
        \end{equation}
        \begin{equation}
            \frac{(s/\tau)^{T}(z/\tau)}{\max\{1, \ -q^{T}x, \ -b^{T}y-h^{T}z \}} < \epsilon_{\text{gap}}
        \end{equation}
        \label{equ3: optiamlity criteria}
    \end{subequations}
    \noindent
    To detect infeasibility, the solver declares primal infeasibility if:
    \begin{subequations}
        \begin{equation}
            -b^{T}y - h^{T}z < \epsilon_{\text{abs}}
        \end{equation}
        \begin{equation}
            \frac{\Vert A^{T}y + G^{T}z \Vert_{2}}{\max\{1, \ \Vert y \Vert_{2} + \Vert z \Vert_{2} \}} < \epsilon_{\text{rel}}
        \end{equation}
        \label{equ3: pinf criteria}
    \end{subequations}
    \noindent
    Similarly, dual infeasibility is declared if:
    \begin{subequations}
        \begin{equation}
            -q^{T}x < \epsilon_{\text{abs}}
        \end{equation}
        \begin{equation}
            \frac{\Vert Qx \Vert_{2}}{\max\{1, \ \Vert x \Vert_{2} \}} < \epsilon_{\text{rel}}
        \end{equation}
        \begin{equation}
            \frac{\Vert Ax \Vert_{2}}{\max\{1, \ \Vert x \Vert_{2} \}} < \epsilon_{\text{rel}}
        \end{equation}
        \begin{equation}
            \frac{\Vert Gx + s \Vert_{2}}{\max\{1, \ \Vert x \Vert_{2} + \Vert s \Vert_{2} \}} < \epsilon_{\text{rel}}
        \end{equation}
        \label{equ3: dinf criteria}
    \end{subequations}
    In this implementation, we set the default values as $\epsilon_{\text{feas}} = \epsilon_{\text{gap}} = \epsilon_{\text{abs}} = \epsilon_{\text{rel}} = 10^{-8}$.

\subsection{Overall Algorithm} \label{subsec: overall algorithm}

    The overall algorithm of the proposed solver is summarized in \textbf{Algorithm \ref{alg: solver algorithm}}. Here, the variable $\texttt{FLAG}$ denotes the solution status indicator, where the $\texttt{enum}$-type status values $\texttt{OPT}$, $\texttt{PINF}$, and $\texttt{DINF}$ represent successful convergence to an optimal solution, primal infeasibility, and dual infeasibility, respectively.

\section{Solver Customization} \label{sec: solver customization}
    
    This section presents the solver customization framework and the code generation tool developed in this study. A custom solver is tailored to handle a specific problem class that shares the same sparsity pattern in the data matrices $Q$, $A$, and $G$, as well as the same underlying structure of the cone $\mathcal{K}$. We refer to such a problem class as a \textit{problem family}, and each individual problem within a problem family as a \textit{problem instance}. The code generation tool allows users to easily obtain a custom solver for a desired problem family. The obtained custom solver is then integrated into the user application to solve problem instances belonging to the same problem family. The overall workflow for the customization and the code generation process is illustrated in Fig. \ref{fig: workflow}.

        \begin{algorithm}[t]
        \caption{Solver Algorithm}
        \begin{algorithmic}[1]
            \REQUIRE Tolerance $\epsilon_{\text{feas}}$, $\epsilon_{\text{gap}}$, $\epsilon_{\text{abs}}$ and $\epsilon_{\text{rel}}$ for termination criteria, and the problem data given as: $Q$, $q$, $A$, $b$, $G$, $h$ and $\mathcal{K}$. \\
            \STATE \textbf{Initialization:} Set $k=0$. Obtain $(x^{(0)}, y^{(0)}, z^{(0)}, s^{(0)})$ from \S \ref{subsec : initialization}. Set $\kappa^{(0)} = 1$, $\tau^{(0)} = 1$.
            \STATE \textbf{loop:}
            \STATE \quad \textbf{Compute residuals:} $(r_{x}, r_{y}, r_{z}, r_{\tau})$ from \eqref{equ3: residuals}.
            \STATE \quad \textbf{Check termination criteria:}
            \STATE \qquad \textbf{if:} \eqref{equ3: optiamlity criteria} holds 
            \STATE \qquad \quad $\texttt{FLAG} \leftarrow \texttt{OPT}$ and \textbf{break}.
            \STATE \qquad \textbf{else if:} \eqref{equ3: pinf criteria} holds 
            \STATE \qquad \quad $\texttt{FLAG} \leftarrow \texttt{PINF}$ and \textbf{break}.
            \STATE \qquad \textbf{else if:} \eqref{equ3: dinf criteria} holds 
            \STATE \qquad \quad $\texttt{FLAG} \leftarrow \texttt{DINF}$ and \textbf{break}.
            \STATE \quad \textbf{Update NT-scaling:} from \eqref{equ3: NNC scaling} and \eqref{equ3: SOC scaling}.
            \STATE \quad \textbf{Factor KKT matrix:} using $LDL^{T}$ factorization.
            \STATE \quad \textbf{Compute search direction:}
            \STATE \qquad Compute affine search direction $\Delta\chi_{a}$.
            \STATE \qquad Obtain centering parameter $\sigma$ from \eqref{equ3: centering parameter}.
            \STATE \qquad Compute nonlinearity error from \eqref{equ3: nonlinearity error}.
            \STATE \qquad Compute combined search direction $\Delta \chi_{c}$.
            \STATE \quad \textbf{Update iterates:} using \eqref{equ3: iterates update}.
            \STATE \quad $k \leftarrow k+1$
            \STATE \textbf{end loop}
            \STATE \textbf{if:} \texttt{FLAG} == \texttt{OPT}
            \STATE \quad \textbf{return} $(x^{\ast}, y^{\ast}, z^{\ast}, s^{\ast}) = (x, y, z, s)^{(k)}/\tau^{(k)}$
            \STATE \textbf{else:}
            \STATE \quad \textbf{return} \texttt{FLAG}
        \end{algorithmic}
        \label{alg: solver algorithm}
    \end{algorithm}

    \subsection{Exploiting Sparsity Patterns of the KKT System}
    
    In this study, the sparsity pattern of the problem family is analyzed offline to perform AMD ordering and \textit{symbolic factorization} which compute the sparsity pattern of the factor $L$ \cite[\S 4]{davis2006direct}. Since these procedures depend solely on the sparsity pattern of $\tilde{K}$, rather than its numerical values, they can be computed once offline and applied to all problem instances within the same problem family. The resulting permutation and sparsity pattern of $L$ are cached in the custom solver code and utilized online during \textit{numerical factorization} which computes the numerical values of the factors \cite[\S 4]{davis2006direct}. This customization framework enables the custom solver to fully rely on static memory allocation of the problem data without the need for dynamic memory allocation (i.e., no use of memory allocation functions such as \texttt{malloc} during implementation). 
    For numerical factorization, the algorithm proposed in \cite{davis2005algorithm} is implemented using nested \texttt{for}-loops, which enables code generation that is nearly independent of problem size \cite{domahidi2012efficient, banjac2017embedded}.

    \subsection{Code Generation and Parsing Information for Updating Problem Instances}

\begin{figure*}[t]
        \centering
        \includegraphics[width=1.0\linewidth]{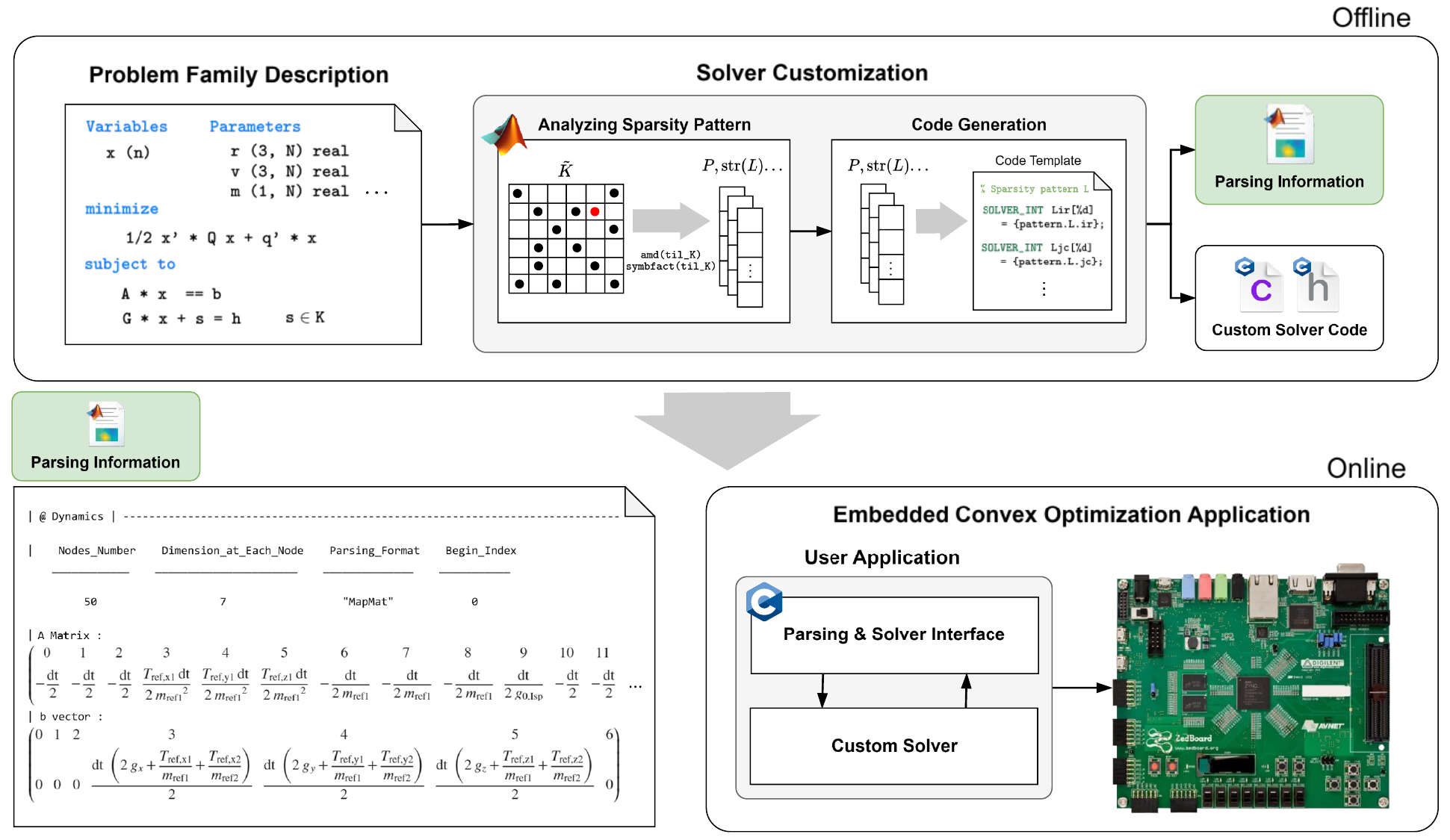}
        \caption{Workflow of the custom solver and code generation tool developed in this study.}
        \label{fig: workflow}
    \end{figure*}

    To enable efficient solver customization, we develop a code generation tool, implemented in MATLAB, that produces a custom solver code tailored to a given problem family. It analyzes the sparsity pattern of the problem and generates C code for a custom solver using a predefined code template. 
    The generated code requires no external dependencies other than the standard library \texttt{math.h} for square root computations, and utilizes only the upper triangular part of the KKT matrix for its linear algebra operations, thereby reducing memory usage compared to using the full matrix.
    The input interface of the code generation tool allows users to provide either the problem data for a specific problem instance or a symbolic description of the problem family using the MATLAB Symbolic Math Toolbox. Through the latter interface, the code generation tool produces not only the custom solver code but also the associated parsing information necessary to update the problem data for varying problem instances. (See bottom-left corner of Fig. \ref{fig: workflow} for an example of the parsing information provided by the code generation tool.)
    This feature particularly facilitates handling of problem data stored in compressed column storage (CCS) format, where direct access to individual elements is nontrivial.
    Furthermore, by providing symbolic expressions for the equations corresponding to each element of the problem data matrices, the tool significantly reduces the user’s workload when using the solver.
    
    It is worth noting that symbolic objects are used solely to represent the equations of each element of problem data and are not involved in any computations required for solver customization. Therefore, the symbolic objects from the MATLAB Symbolic Math Toolbox can, in principle, be replaced with other objects such as plain strings. In this study, the toolbox is adopted purely for implementation convenience.

\section{Benchmark Tests} \label{sec: benchmark test}
    
    Benchmark tests were conducted to evaluate the performance of the developed solver across a wide range of convex optimization problems, including the Maros-Mészáros QP problem set, portfolio optimization, and LASSO. For comparison, we used the following convex conic optimization solvers: the open-source IPM solver ECOS \cite{domahidi2013ECOS}, the commercial IPM solver MOSEK \cite{mosek2024}, and the open-source ADMM-based solver SCS \cite{odonoghue2021operator}.
    All solvers were benchmarked using a default termination tolerance of $10^{-8}$. For a fair comparison, a relaxed tolerance of $10^{-4}$ was also applied to SCS, considering the poor tail convergence characteristics of the ADMM algorithm \cite{he2012on}. Pre-solving and pre-conditioning were disabled to ensure that all solvers address the same problem formulation.   
    All tests were conducted in single-threaded mode on a desktop computer equipped with a 13th Gen Intel Core i5-13600K CPU running at 3.5GHz and 64GB of RAM.

    In the benchmark tests, the performance of the solvers is evaluated using the \textit{performance profile} \cite{dolan2002benchmarking}, which provides a comprehensive comparison of solver performance across diverse problem collections. Let us consider the set of problem collections $\mathbf{P}$ and the set of solvers $\mathbf{S}$ for the current benchmark test. The \textit{relative performance ratio} for solver $s \in \mathbf{S}$ and problem $p \in \mathbf{P}$ is defined as,
    \begin{equation}
        r_{p, s} = \frac{t_{p, s}}{\displaystyle \min_{s \in \mathbf{S}} \left\{ t_{p, s} \right\} }
    \end{equation}
    \noindent
    where $t_{p, s}$ is the time taken by solver $s$ to solve problem $p$. The relative performance profile for solver $s$ is obtained by plotting the function $\rho_{s}^{r}$ with respect to the relative performance ratio $r_{p, s}$ for $s$. The function $\rho_{s}^{r}$ is defined as:
    \begin{equation}
        \rho^{r}_{s}(\tau) = \frac{1}{N_{\mathbf{P}}} \sum_{p \in \mathbf{P}} \mathcal{I}(r_{p, s} ; \tau)
    \end{equation}
    \noindent
    where $N_{\mathbf{P}}$ represents the number of problems in the problem set $\mathbf{P}$, and the indicator function $\mathcal{I}$ is defined as:
    \begin{equation}
        \mathcal{I}(x;\tau) = \begin{cases} 1 & \text{if } x \leq \tau \\ 0& \text{otherwise} \end{cases}
    \end{equation}
    This performance metric provides insights into both the relative performance of each solver with respect to the fastest solver and its robustness across the problem set.

     We further employ \textit{absolute performance profile} \cite{goulart2024clarabel} to evaluate the absolute performance of the solver, independent of the composition of the solver set $\mathcal{S}$. The absolute performance profile of the solver $s$ is a plot of the function $\rho_{s}^{a}$ with respect to the computation time $t_{p, s}$, where the function $\rho_{s}^{a}$ is defined as:

    \begin{equation}
        \rho_{s}^{a}(\tau) = \frac{1}{N_{\mathbf{P}}} \sum_{p \in \mathbf{P}} \mathcal{I}(t_{p, s} ; \tau)
    \end{equation}

    \subsection{Maros-Mészáros QP Problem Set}
    
    Maros-Mészáros QP problem set \cite{maros1999repository} contains 138 QP problems with wide range of sizes and conditions. In particular, it includes many challenging cases arising from numerical ill-conditioning or rank deficiency. In this benchmark test, the results for SCS solver with full tolerance level are excluded due to its excessively long computation times compared to other solvers, making the comparison meaningless. The prolonged computation times are attributed to the poor tail convergence characteristics of ADMM.

\begin{figure*}
        \centering
        \subfloat[Relative Performance Profile]
        {\includegraphics[width = 0.49\textwidth]{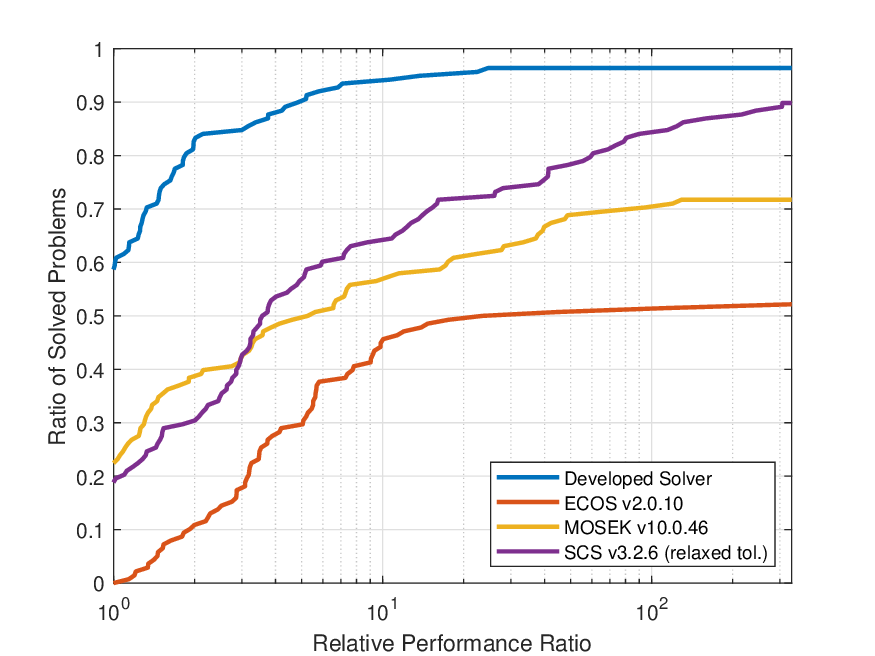}}%
        \label{subfig: Maros-Meszaros relative}
        \hfil
        \subfloat[Absolute Performance Profile]
        {\includegraphics[width = 0.49 \textwidth]{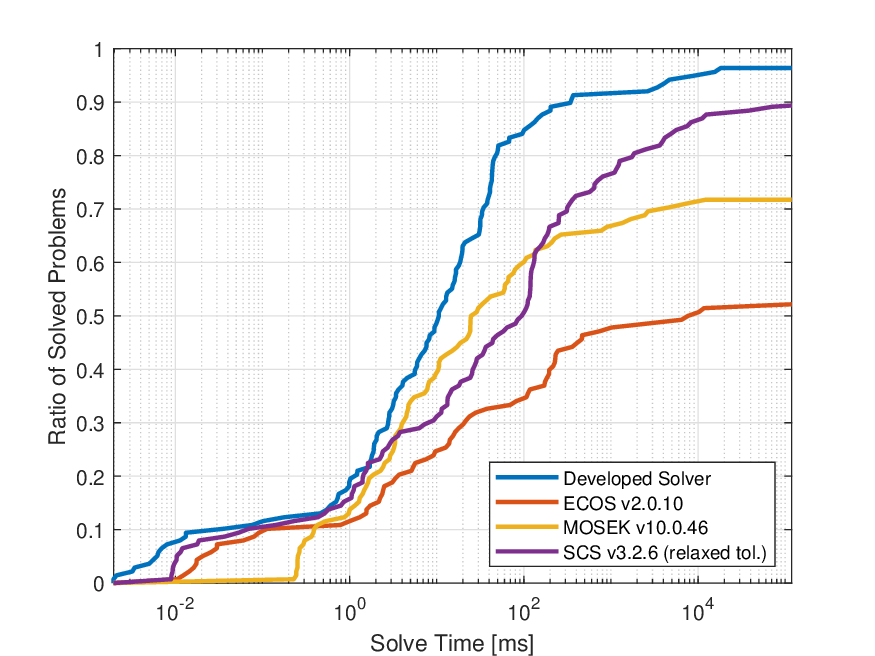}}%
        \label{subfig: Maros-Meszaros absolute}
        \caption{Performance profiles for Maros-Mészáros QP problem set}
        \label{fig: Maros-Meszaros}
    \end{figure*}

    \textit{Results}—Figure \ref{fig: Maros-Meszaros} presents the benchmark results for the Maros-Mészáros QP problem set. The results show that the solver developed in this study outperforms other solvers in terms of both computation time and success rate across most problem cases. In the case of ECOS and MOSEK, the low performance and high failure rate are attributed to their algorithms being limited to handling only linear cost functions, which necessitates the reformulation of a quadratic objective term into a conic constraint as in \eqref{equ3: linear conic program}. In particular, this requirement leads to numerical failure in a substantial number of cases when computing $Q^{1/2}$, due to the ill-conditioning of the quadratic cost matrix $Q$. SCS exhibits a low failure rate due to its ability to directly handle quadratic cost functions. However, it suffers from slow convergence, resulting in degraded performance even under the relaxed tolerance.

    \begin{figure*}
        \centering
        \subfloat[Relative Performance Profile]
        {\includegraphics[width = 0.49\textwidth]{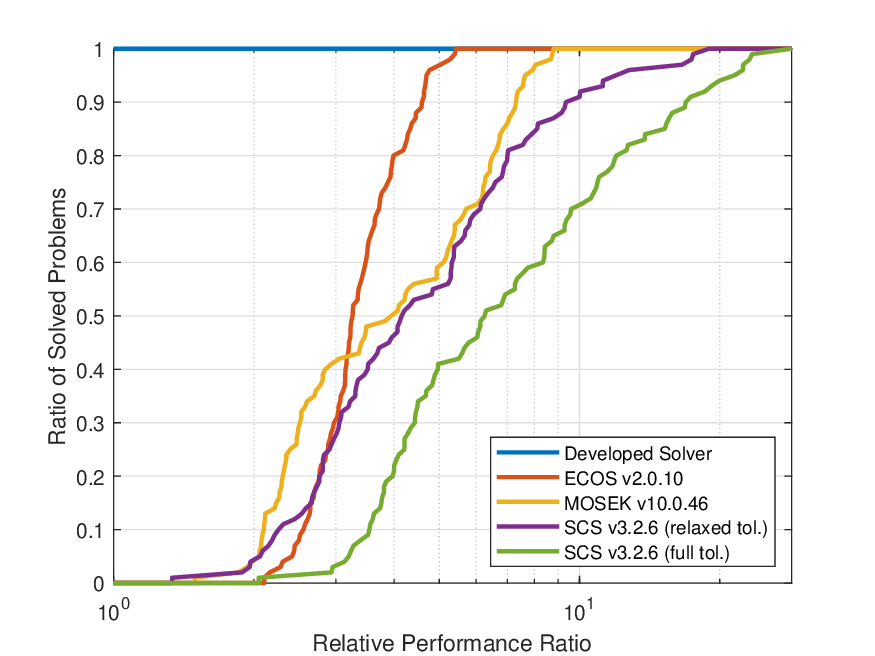}}%
        \label{subfig: Portfolio relative}
        \hfil
        \subfloat[Absolute Performance Profile]
        {\includegraphics[width = 0.49 \textwidth]{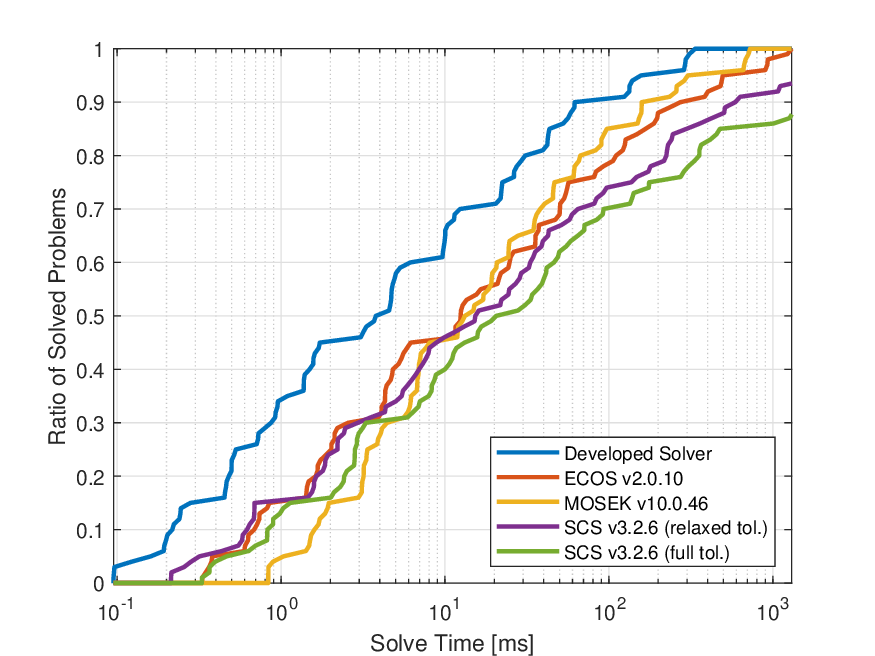}}%
        \label{subfig: Portfolio absolute}
        \caption{Performance profiles for portfolio optimization problem}
        \label{fig: Portfolio}
    \end{figure*}

    \subsection{Portfolio Optimization}

    Portfolio optimization problem seeks to find an asset allocation that maximizes the risk-adjusted return \cite[\S4.7.6]{boyd2004convex}.

    \begin{equation}
        \begin{split}
            \underset{\omega}{\text{maximize}} & \quad \bar{\mu}^{T}\omega - \rho \left( \omega^{T}\Sigma \omega \right) \\
            \text{subject to} & \quad \mathbf{1}^{T}\omega = 1 \\
            & \quad \omega \geq 0
        \end{split}
    \end{equation}
    where $\omega \in \mathbb{R}^{n}$ represents the portfolio, $\bar{\mu} \in \mathbb{R}^{n}$ is the vector of expected return, $\rho > 0$ is the risk-aversion parameter, and $\Sigma \in \mathbb{S}^{n}_{+}$ represents the risk model covariance matrix. This risk covariance model is typically expressed in factor form as:
    \begin{equation}
        \Sigma = FF^{T} + \Pi
    \end{equation}
    \noindent
    where $F \in \mathbb{R}^{n \times k}$ ($k < n$) is the factor-loading matrix, and $\Pi \in \mathbb{R}^{n \times n}$ is a diagonal matrix representing idiosyncratic risk. To simplify the problem, we introduce a new variable $\zeta = F^{T}\omega$, which transforms the original formulation into:
    \begin{equation}
        \begin{split}
            \underset{\omega, \zeta}{\text{maximize}} & \quad \frac{1}{2}\omega^{T}\Pi\omega + \frac{1}{2} \zeta^{T}\zeta -\frac{1}{2\rho} \bar{\mu}^{T}\omega \\
            \text{subject to} & \quad \zeta = F^{T}\omega, \quad \mathbf{1}^{T}\omega = 1  \\
            & \quad \omega \geq 0
        \end{split} 
    \end{equation} 
    making the Hessian of the quadratic cost function diagonal. The problem instances for benchmark test are generated following the similar approach used in \cite{stellato2020osqp}. The instances are generated by increasing number of factors $k = [5, \ 10, \ 20, \ 50, \ 100]$ and ratios of $n/k = [10, \ 20, \ 50, \ 100]$. For each problem size, five different problems are generated, thus total 100 instances are used for the benchmark test. The elements of the $F$ matrix are sampled as  $F_{i, j} \sim \mathcal{N}(0, 1)$ with 70\% nonzero elements, the elements of $\Pi$ matrix are sampled as $\Pi_{i, i} \sim \mathcal{U}[0, \sqrt{k}]$, and mean return $\bar{\mu}$ is sampled as $\bar{\mu}_{i} \sim \mathcal{N}(0, 1)$. The total number of nonzero components in the problem data matrices ranges from 338 to 54,355.

    \textit{Results}—Figure \ref{fig: Portfolio} shows the benchmark results for portfolio optimization problems. We can observe that the solver developed in this study outperforms the other solvers across all problem instances. 
    The problem instances in this benchmark set correspond to small- to medium-scale problems, where IPM-based solvers generally exhibit better performance than the ADMM-based solver. Moreover, among the IPM-based solvers, the developed solver achieves higher computational efficiency owing to its ability to directly handle quadratic objective functions. This capability allows matrix factorization to be performed on a smaller and sparser KKT matrix than those of the other IPM-based solvers, ECOS and MOSEK, which are limited to linear objective functions. Therefore, these results demonstrate the effectiveness of the algorithmic formulation for handling quadratic cost functions, as well as the strength of IPMs in small- to medium-scale problems.

    \subsection{LASSO}

    \begin{figure*}
        \centering
        \subfloat[Relative Performance Profile]
        {\includegraphics[width = 0.49\textwidth]{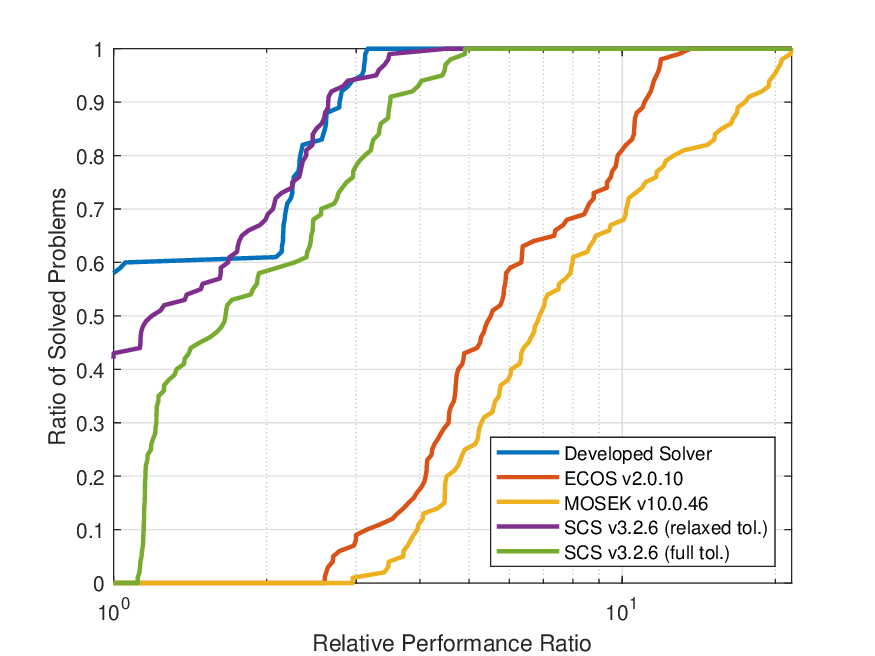}}%
        \label{subfig: LASSO relative}
        \hfil
        \subfloat[Absolute Performance Profile]
        {\includegraphics[width = 0.49 \textwidth]{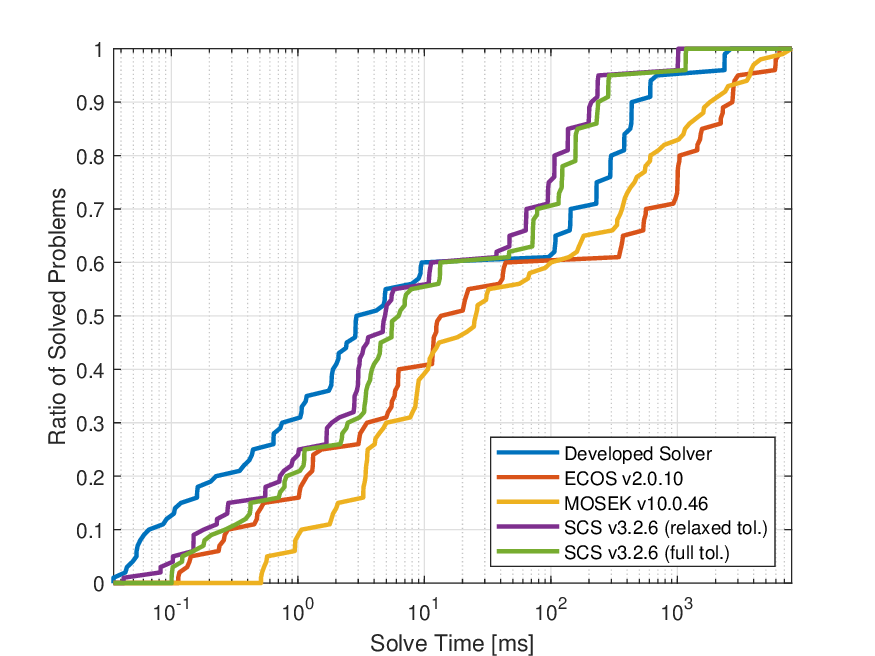}}%
        \label{subfig: LASSO absolute}
        \caption{Performance profiles for LASSO problem}
        \label{fig: LASSO}
    \end{figure*}

    The least absolute shrinkage and selection operator (LASSO) is a technique for sparse linear regression with an $\mathcal{L}_{1}$ regularization term in the cost function \cite{tibshirani1996regression, candes2008enhancing}.
    \begin{equation}
        \begin{split}
            \underset{\theta}{\text{minimize}} & \qquad \frac{1}{2}\Vert X \theta - \bar{y} \Vert_{2}^{2} + \bar{\lambda}\Vert \theta \Vert_{1}
        \end{split}
    \end{equation}
    where $\theta\in\mathbb{R}^{n}$ is the vector of parameters, $X \in \mathbb{R}^{m\times n}$ is the data matrix, and $\bar{\lambda} > 0$ is the weighting parameter. For LASSO, the problem instances for benchmark test are also generated following the method used in \cite{stellato2020osqp}. We increase the number of regression parameter $n = [10, \ 50, \ 100, \ 500, \ 700]$ and the ratio $m/n = [2, \ 5, \ 10, 20]$. For each problem size, five different problems are generated, thus total 100 instances are used for the benchmark test. The elements of matrix $X$ are sampled as $X_{i, j} \sim \mathcal{N}(0, 1)$ with 15\% nonzero elements. To construct the vector $\bar{y}$, we construct the true sparse vector $v \in \mathbb{R}^{n}$ to be learned as:
    \begin{equation}
        v_{i} \sim \begin{cases} 0 & \text{sampled with probability of 0.5} \\ \mathcal{N}(0, 1/n) & \text{otherwise} \end{cases}
    \end{equation}
    Then, the vector $\bar{y}$ can be constructed with noise $\epsilon_{i} \sim \mathcal{N}(0, 1)$ as follows:
    \begin{equation}
        \bar{y}_{i} = \sum_{j=1}^{n}X_{i, j}v_{j} + \epsilon_{i} \quad (i=1, 2, \cdots, m)
    \end{equation}
    Finally, the weighting parameter $\bar{\lambda}$ is chosen as $(1/5)\Vert X^{T}\bar{y} \Vert_{\infty}$. The total number of nonzero components in the problem data matrices ranges from 108 to 1,395,786.

    \textit{Results}—Figure \ref{fig: LASSO} presents the benchmark results for LASSO problems. The solver developed in this study outperforms other solvers on most instances, particularly those of small- to medium-scale. However, for large-scale problems—specifically those with $n > 500$, $m > 1000$ and more than 73,628 nonzero elements—the ADMM-based solver SCS shows better performance.
    This is primarily because the computational cost for matrix factorization significantly increases with problem size and sparsity. 
    Since ADMM-based solvers requires matrix factorization only once at initialization, they can achieve favorable performance on large-scale problems due to their low computational cost per iteration.
    Nevertheless, the solver developed in this study still demonstrates better performance than the other IPM-based solvers attributed to its more efficient algorithmic formulation.

\section{Numerical Experiments on Embedded System} \label{sec: embedded system experiments}
    The solver developed in this study is specifically designed for embedded real-time optimization. Therefore, we evaluated its performance on G\&C problems in the embedded system environment. The evaluations were conducted on an AMD Zynq-7000 system-on-chip featuring an ARM Cortex-A9 processor running at a maximum clock frequency of 650MHz with 512MB of RAM. For comparison, we used ECOS and SCS, as both are implemented in C and distributed as open-source software, facilitating deployment on embedded systems. All solvers were compiled with the \texttt{-O3} optimization level to maximize execution efficiency.

    \subsection{Trajectory Optimization for Mars Landing}
    
    The first numerical experiment considers a trajectory optimization problem for Mars landing. This problem is formulated as the following SOCP problem by applying convexification techniques, including lossless convexification and change of variables \cite{acikmese2007convex}.
    \begin{equation}
        \begin{split}
            \underset{x, z, u, \sigma}{\text{minimize}} & \quad -z_{N}\\
            \text{subject to} & \quad x_{k+1} = A_{k}x_{k} + B_{k}u_{k} + g \  \forall k \in \{ 0, 1, \cdots, N-1\}\\
            & \quad z_{k+1} = z_{k} - \alpha \sigma_{k} \Delta t \quad \forall k \in \{ 0, 1, \cdots, N-1\} \\
            & \quad x_{0} = x_{\text{init}}, \ z_{0} = \ln(m_{\text{wet}}), \ x_{N} = {\bf{0}}_{6\times1} \\
            & \quad \mu_{1, k} \left[ 1 - (z_{k} - z_{0, k}) \right] \leq \sigma_{k} \quad \forall k \in \{ 0, 1, \cdots, N\} \\
            & \quad \sigma_{k} \leq \mu_{2, k} \left[ 1 - (z_{k} - z_{0, k}) \right] \quad \forall k \in \{ 0, 1, \cdots, N\} \\
            & \quad \Vert u_{k} \Vert \leq \sigma_{k} \quad \forall k \in \{ 0, 1, \cdots, N\} \\
            & \quad e_{3}^{T}u_{k} \geq \sigma_{k}\cos{\theta_{\max}} \quad \forall k \in \{ 0, 1, \cdots, N\} \\
            & \quad z_{N} \geq \ln(m_{\text{dry}})
        \end{split}
    \end{equation}
    where $x_{k} \in \mathbb{R}^{6}$ is the position and velocity of the lander, $z_{k} \in \mathbb{R}$ is the natural logarithm of the lander's mass, $u_{k} \in \mathbb{R}^{3}$ is the acceleration generated by the lander's thrusters, and  $\sigma_{k} \in \mathbb{R}$ is the slack variable at time node $k$. The parameters $m_{\text{dry}}$ and $m_{\text{wet}}$ represent the dry mass and wet mass of the lander, respectively, $e_{3} \in \mathbb{R}^{3}$ is the $z-$axis basis vector of the landing inertial coordinate frame, $\theta_{\max}$ is the maximum tilt angle of the lander from $e_{3}$ vector, and $x_{\text{init}} \in \mathbb{R}^{6}$ is the initial states of the lander. We define $z_{0, k} = \ln(m_{\text{wet}} - \alpha \rho_{2}k\Delta t)$, $\mu_{1, k} = \rho_{1}e^{-z_{0, k}}$, $\mu_{2, k} = \rho_{2} e^{-z_{0, k}}$ where $\rho_{1}$ and $\rho_{2}$ is the minimum and maximum thrust magnitudes, respectively, $\alpha$ is the fuel depletion parameter, and $\Delta t$ is the time step between nodes. The matrices $A_{k}$ and $B_{k}$, and the vector $g$ in the dynamics are given by:
    \begin{equation}
        \begin{split}
            A_{k} = \begin{bmatrix} I_{3} & I_{3}\Delta t \\ {\bf{0}}_{3\times 3} & I_{3} \end{bmatrix}, \quad B_{k} = \begin{bmatrix} (\frac{1}{2}\Delta t^{2})I_{3} \\ I_{3} \Delta t \end{bmatrix} \\
            g = \begin{bmatrix} 0 & 0 & -\frac{1}{2}g_{0}\Delta t^{2} & 0 & 0 & -g_{0}\Delta t \end{bmatrix}^{T}
        \end{split}
    \end{equation}
    where $I_{3}$ denotes 3$\times$3 identity matrix, $g_{0}$ represents the gravitational acceleration on the Martian surface, and the following parameter values were used in this numerical experiment.
    \begin{equation*}
        \centering
        \begin{split}
            x_{\text{init}} = [ 200.0, \ 0.0, \ 800.0, \ -35.0, \ 0.0, \ -75.0]^{T} \\
            \rho_{1} = 7440, \ \rho_{2} = 18600, \ g_{0} = 3.7114, \ \theta_{\max} = \pi/12 \\
            m_{\text{dry}} = 1505, \ m_{\text{wet}} = 1905, \ \alpha = 4.53 \times 10^{-4}
        \end{split}
    \end{equation*}
    
    Moreover, to validate the infeasibility detection capability of the developed solver, numerical experiments were conducted with both a feasible and an infeasible case. Since the Mars landing problem becomes infeasible if the flight time $t_{f}$ is shorter than the minimum admissible value \cite{acikmese2007convex, dueri2017customized}, we selected $t_{f} = 48$s as the feasible case and $t_{f} = 25$s as the infeasible case. The feasibility of these cases was confirmed by examining the solutions returned by ECOS.
    For both cases, ten problem instances were generated by increasing the number of nodes $N$ from 25 to 500.
    We set the termination tolerance to $10^{-8}$ for the developed customized solver and ECOS, while a relaxed tolerance of $10^{-4}$ was used for SCS, considering the convergence characteristics of ADMM.

    \begin{table}[t]
        \caption{Differences in Solutions for Feasible Mars Landing Problems}
        \centering
        \renewcommand{\arraystretch}{1.5}
        \resizebox{\linewidth}{!}{
        \textcolor{black}{
        \begin{tabular}{ccccc}
            \hline\hline
             Differences & Min. & Max. & Mean & Median  \\
             \hline
             $\Vert x^{\ast} - x^{\ast}_{\text{ECOS}} \Vert_{\infty}$ & 7.899e-07 & 1.126e-01 & 5.325e-02 & 6.759e-02 \\
             $\Vert x^{\ast} - x^{\ast}_{\text{SCS}} \Vert_{\infty}$  & 2.848e+00 & 1.344e+01 & 8.470e+00 & 8.783e+00 \\
             $\Vert x^{\ast}_{\text{ECOS}} - x^{\ast}_{\text{SCS}} \Vert_{\infty}$  & 2.925e+00 & 1.343e+01 & 8.510e+00 & 8.819e+00 \\
            \hline\hline
        \end{tabular}
        } }
        \label{tab: sol_diff_lcvx}
    \end{table}

    \begin{table}[t]
        \caption{Residuals of the Developed Solver's Solution\\ for Feasible Mars Landing Problems}
        \centering
        \renewcommand{\arraystretch}{1.5}
        \resizebox{\linewidth}{!}{
        \textcolor{black}{
        \begin{tabular}{ccccc}
            \hline\hline
             Residuals & Min. & Max. & Mean & Median  \\
             \hline
             $\Vert Ax^{\ast}-b\Vert_{\infty}$ & 5.733e-12 & 3.621e-11 & 1.734e-11 & 1.425e-11 \\
             $\Vert Gx^{\ast} + s^{\ast} - h \Vert_{\infty}$  & 4.489e-11 & 7.645e-05 & 1.457e-05 & 5.616e-08 \\
             $s^{\ast T}z^{\ast}$  & 5.228e-10 & 2.862e-08 & 1.655e-08 & 1.825e-08 \\
            \hline\hline
        \end{tabular}
        } }
        \label{tab: resd_lcvx}
    \end{table}

    \textit{Results}—Table \ref{tab: sol_diff_lcvx} presents the differences in the converged solutions among the solvers for the feasible problem instances, where $x^{\ast}$, $x^{\ast}_{\text{ECOS}}$, and $x^{\ast}_{\text{SCS}}$ denote the solutions obtained by the developed solver, ECOS, and SCS, respectively. When compared with ECOS, the developed solver produced solutions that differed by values on the order of $10^{-7}$ to $10^{-1}$. These small discrepancies indicate that the developed solver produces accurate solutions consistent with those obtained by the existing solver. On the other hand, relatively larger discrepancies on the order of $10^{0}$-$10^{1}$ were observed for SCS, since it used a relaxed termination tolerance to accommodate its convergence characteristics.

    The solution accuracy is further assessed through the residual values for feasibility and optimality as shown in Table \ref{tab: resd_lcvx}. Specifically, $\Vert Ax^{\ast}-b\Vert_{\infty}$ and $\Vert Gx^{\ast} + s^{\ast} - h\Vert_{\infty}$ quantify the maximum feasibility residuals for the equality and inequality constraints, respectively, while $s^{\ast T}z^{\ast}$ represents the duality gap, serving as a measure of optimality. Here, $A$, $b$, and $G$, $h$ denote the matrices and vectors corresponding to the equality and inequality constraints, respectively, which appear when the problem is expressed in the form of \ref{equ2: primal problem}. The results show that the residual values on the order of the prescribed tolerance verify the accuracy of the solutions obtained by the developed solver.

\begin{figure}[t]
        \centering
         \includegraphics[width=\linewidth]{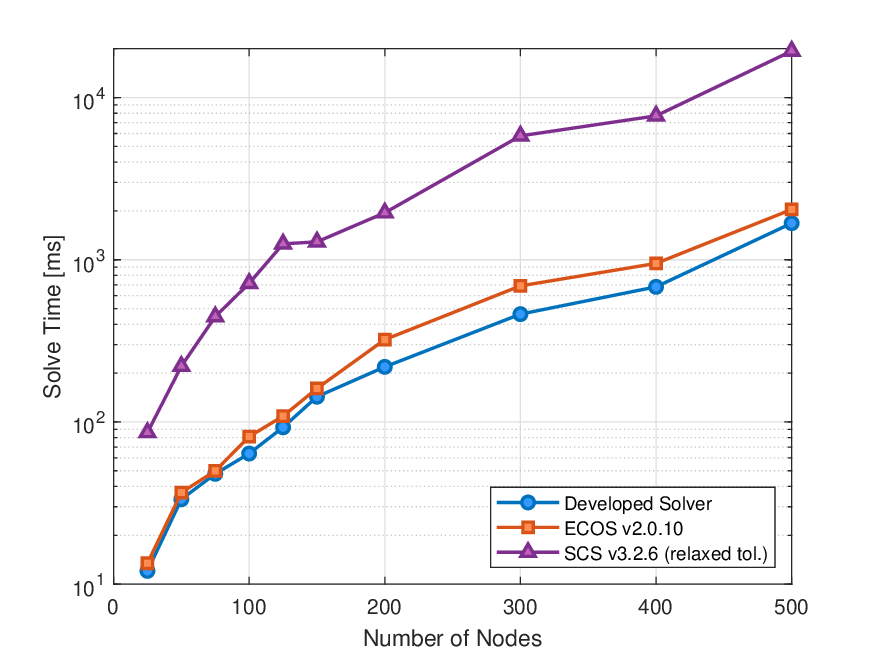}
        \caption{Computation time results for the feasible Mars landing trajectory optimization examples}
        \label{fig: LCVX emb}
    \end{figure}

    Figure \ref{fig: LCVX emb} presents the computation time results for the feasible Mars landing trajectory optimization problems. It can be observed that the developed solver and ECOS achieve faster onboard computation times than SCS even under tighter termination tolerances, demonstrating the effectiveness of IPMs for G\&C applications. For small-scale instances ($N \leq 150$), ECOS and the developed solver exhibit comparable performance. However, for larger-scale instances, the solver developed in this study demonstrates faster performance. This improvement is attributed to the elimination of unnecessary logic and the enhanced efficiency of data memory management achieved through solver customization.

    Lastly, Table \ref{tab: infeas_LCVX} presents the results demonstrating the infeasibility detection capability of the developed solver. It shows the converged values of $\kappa^{\ast}$ and $\tau^{\ast}$ together with the number of solver iterations for the infeasible instances. As described in Section \ref{subsec : homogeneous embedding}, the condition $\kappa^{\ast} > 0$ with $\tau^{\ast} = 0$ provides a certificate of infeasibility for the given problem. The resulting values of $\kappa^{\ast}$ and $\tau^{\ast}$, which satisfy the infeasibility conditions within the numerical tolerance, confirm that the solver successfully detects infeasibility. Moreover, it is observed that the solver achieves convergence to a certificate of infeasibility in only 7–8 iterations. This implies that, even for infeasible instances, the algorithm converges within a finite number of iterations and explicitly provides infeasibility information, which can be exploited to establish contingency strategies in G\&C applications under abnormal conditions.
    
    \begin{table}[t]
        \caption{Infeasibility Detection Results for Infeasible \\ Mars Landing Problems}
        \centering
        \renewcommand{\arraystretch}{1.4}
        \begin{tabularx}{0.88\linewidth}{YYYYY}
            \hline\hline
            $N$ & $\kappa^{\ast}$ & $\tau^{\ast}$ & $\kappa^{\ast}/\tau^{\ast}$ & Iteration \\
            \hline
            25  & 1.742e-01 & 1.131e-08 & 1.541e+07 & 8 \\
            50  & 1.760e-01 & 1.168e-08 & 1.506e+07 & 8 \\
            75  & 1.742e-01 & 1.131e-08 & 1.732e+05 & 8 \\
            100 & 1.765e-01 & 9.943e-07 & 1.775e+05 & 7 \\
            125 & 1.762e-01 & 9.752e-07 & 1.807e+05 & 7 \\
            150 & 1.749e-01 & 9.090e-07 & 1.924e+05 & 7 \\
            200 & 1.709e-01 & 8.097e-07 & 2.111e+05 & 7 \\
            300 & 1.665e-01 & 7.352e-07 & 2.264e+05 & 7 \\
            400 & 1.624e-01 & 6.645e-07 & 2.444e+05 & 7 \\
            500 & 1.690e-01 & 8.418e-09 & 2.007e+07 & 7 \\
            \hline\hline
        \end{tabularx}
        \label{tab: infeas_LCVX}
    \end{table}

    \subsection{Quadcopter Model Predictive Control}
    
    Finally, we consider an MPC problem for quadcopter control. This problem is formulated as a finite-horizon constrained optimal control problem, which can be expressed as the following QP.
    
    \begin{equation}
        \begin{split}
            \underset{x, u}{\text{minimize}} & \quad (x_{N} - x_{r})^{T} Q_{N} (x_{N} - x_{r}) \  \\
            & \quad +\sum_{k=0}^{N-1} \left[ (x_{k} - x_{r})^{T}Q(x_{k} - x_{r}) + u_{k}^{T}Ru_{k} \right] \\
            \text{subject to} & \quad x_{k+1} = A_{k}x_{k} + B_{k}u_{k} \quad \forall k \in \{ 0, 1, \cdots, N-1\} \\
            & \quad x_{0} = x_{\text{init}} \\
            & \quad x_{\min} \leq x_{k} \leq x_{\max} \quad \forall k \in \{ 0, 1, \cdots, N\} \\
            & \quad u_{\min} \leq u_{k} \leq u_{\max} \quad \forall k \in \{ 0, 1, \cdots, N-1\}
        \end{split}
    \end{equation}
    where $u\in\mathbb{R}^{4}$ denotes the control input, and $x \in \mathbb{R}^{12}$ denotes the state of the quadcopter, consisting of its position, velocity, attitude angles, and angular rates. The matrix $Q_{N} \in \mathbb{S}_{+}^{12}$ defines the final stage cost, and the matrices $Q \in \mathbb{S}_{+}^{12}$ and $R \in \mathbb{S}_{++}^{4}$ define the stage cost for the states and control inputs, respectively. The parameter $x_{\text{init}} \in \mathbb{R}^{12}$ represents the initial state, and $x_{r} \in \mathbb{R}^{12}$ is the reference state, while $x_{\min}, x_{\max} \in \mathbb{R}^{12}$ and $u_{\min}, u_{\max} \in \mathbb{R}^{4}$ represent the state and control constraints, respectively. We used the quadcopter dynamics model from the MPC benchmark collection presented in \cite{kouzoupis2015towards}. Ten problem instances were generated by increasing the number of nodes $N$ from 15 to 100, and the following parameter values were used in the numerical experiment.
    \begin{equation*}
        \begin{split}
            Q = Q_{N} &= \text{diag}(10, 10, 10, 5, 5, 5, 0, 0, 10, {\bf{0}}_{1\times 3}) \\ 
            R &= 0.1 \times I_{4} \\
            x_{\text{init}} &= \begin{bmatrix} 1.0 & 0.0 & 0.5 & {\bf{0}}_{1\times 9} \end{bmatrix}^{T}\\
            \quad x_{r} &= \begin{bmatrix} {\bf{0}}_{1\times 6} & 0.0 & 0.0 & 1.0 & {\bf{0}}_{1\times 3} \end{bmatrix}^{T} \\
        \end{split}
    \end{equation*}
    \begin{equation*}
        \begin{split}
            x_{\min, i} = \begin{cases} -1.0 & i = 3 \\ -\pi/6 & i = 7, 8 \\ -\infty & \text{otherwise} \end{cases},& \ x_{\max, i} = \begin{cases} -\pi/6 & i = 7, 8 \\ \infty & \text{otherwise} \end{cases}\\
            u_{\min} = (-0.9916) {\bf{1}}_{4\times 1},& \ u_{\max} = (2.4084){\bf{1}}_{4\times 1}
        \end{split}
    \end{equation*}
    We set the termination tolerance to $10^{-6}$ for the developed customized solver and ECOS, while a relaxed tolerance of $10^{-4}$ was used for SCS.

    \textit{Results}—Table \ref{tab: sol_diff_quadmpc} presents the solution differences among the solvers for the quadcopter MPC problems. As observed in the results, the solution obtained by the developed solver differed from those of existing solvers by the order of $10^{-3}$-$10^{-5}$, validating its accuracy. This accuracy is further corroborated by the residual values which remain at the prescribed tolerance level, as shown in Table \ref{tab: resd_mpc}.

\begin{figure}[t]
        \centering
        \includegraphics[width=\linewidth]{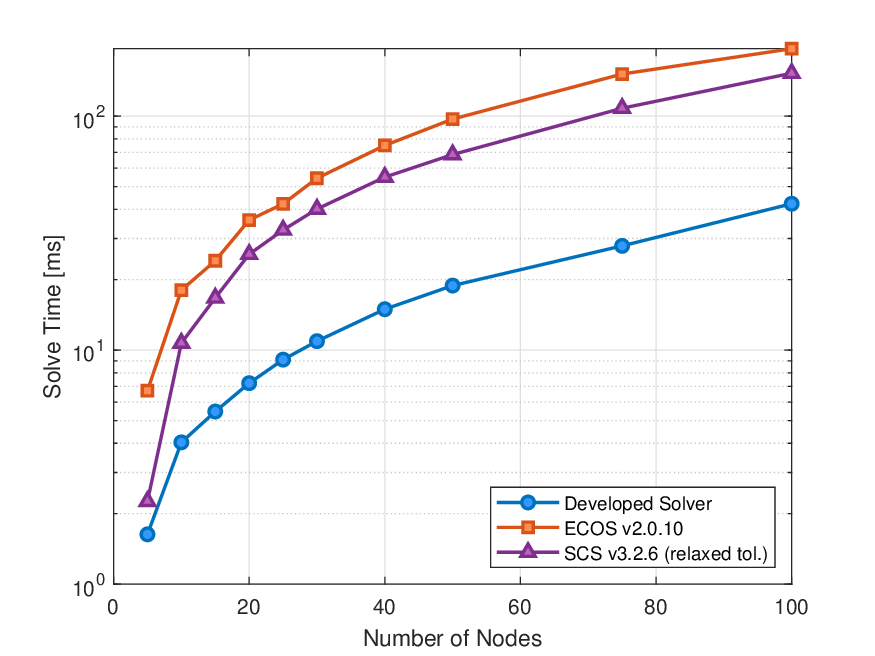}
        \caption{Computation time results for the quadcopter MPC examples}
        \label{fig: MPC emb}
    \end{figure}

    Figure \ref{fig: MPC emb} shows the computation time results for the quadcopter MPC problems. The results demonstrate that the solver developed in this study outperforms the other solvers. Such performance is mainly attributed to the efficient algorithmic formulation that directly handles quadratic cost functions and the fast convergence properties of IPMs. This capability is particularly important in MPC applications, where the optimization problem should be solved repeatedly in real-time based on feedback from the updated states of the quadcopter. 

    \begin{table}[t]
        \caption{Differences in Solutions for Quadcopter MPC Problems}
        \centering
        \renewcommand{\arraystretch}{1.5}
        \resizebox{\linewidth}{!}{
        \begin{tabular}{ccccc}
            \hline\hline
             Differences & Min. & Max. & Mean & Median  \\
             \hline
             $\Vert x^{\ast} - x^{\ast}_{\text{ECOS}} \Vert_{\infty}$ & 6.152e-05 & 4.093e-03 & 1.007e-03 & 6.320e-04 \\
             $\Vert x^{\ast} - x^{\ast}_{\text{SCS}} \Vert_{\infty}$  & 1.983e-04 & 4.848e-03 & 1.295e-03 & 9.675e-04 \\
             $\Vert x^{\ast}_{\text{ECOS}} - x^{\ast}_{\text{SCS}} \Vert_{\infty}$  & 1.091e-04 & 1.805e-03 & 1.054e-03 & 1.076e-03 \\
            \hline\hline
        \end{tabular}
        }
        \label{tab: sol_diff_quadmpc}
    \end{table}

    \begin{table}[t]
        \caption{Residuals of the Developed Solver's Solution \\ for Quadcopter MPC Problems}
        \centering
        \renewcommand{\arraystretch}{1.5}
        \resizebox{\linewidth}{!}{
        \begin{tabular}{ccccc}
            \hline\hline
             Residuals & Min. & Max. & Mean & Median  \\
             \hline
             $\Vert Ax^{\ast}-b\Vert_{\infty}$ & 5.965e-09 & 1.079e-07 & 2.133e-08 & 6.698e-09 \\
             $\Vert Gx^{\ast} + s^{\ast} - h \Vert_{\infty}$  & 8.882e-16 & 2.282e-09 & 9.309e-10 & 6.729e-10 \\
             $s^{\ast T}z^{\ast}$  & 7.540e-06 & 9.970e-05 & 2.559e-05 & 1.304e-05 \\
            \hline\hline
        \end{tabular}
        }
        \label{tab: resd_mpc}
    \end{table}
 
\section{Conclusion} \label{sec: conclusion}
    In this work, we have developed a customized IPM-based solver specifically tailored for embedded real-time convex optimization. The developed solver is capable of handling problems within the class of SOCPs, while supporting infeasibility detection through a homogeneous embedding formulation. Unlike conventional homogeneous self-dual embedding approaches, which are limited to problems with linear cost functions, we adopt a homogeneous embedding framework for MCPs from previous work which can directly accommodate quadratic objective terms. This formulation enables efficient and robust algorithm design by avoiding the reformulation of quadratic objectives into conic constraints.
    To support embedded system deployment, the solver is implemented in C with complete static memory allocation, making it suitable for resource-constrained systems. A code generation tool is also developed to automate the solver customization process for a given problem family, exploiting structural sparsity offline and minimizing computation overhead during online execution.
    Numerical benchmark results demonstrate that the developed solver outperforms existing IPM- and ADMM-based solvers on small- to medium-scale problems in terms of computation time and reliability.
    Although the developed solver shows performance degradation for some large-scale problems due to the limited scalability inherent in IPM algorithm, it still achieves better performance than the other IPM-based solvers due to its efficient algorithmic design.
    Moreover, the results on embedded hardware platform experiments highlight its effectiveness for onboard optimization in G\&C applications such as MPC and trajectory optimization.

\section*{Acknowledgments}
This work was partly supported by Theater Defense Research Center funded by Defense Acquisition Program Administration under Grant UD240002SD, and by the Korea Aerospace Research Institute (KARI)’s own research project titled “Development of Control Performance Analysis Program for Spaceplanes”. The authors acknowledge the use of ChatGPT to improve the syntax and grammar of several paragraphs in the manuscript. The manuscript was originally drafted by the authors, with the AI system employed to enhance clarity and readability without modifying the technical content.

\bibliographystyle{IEEEtran}
\bibliography{ref}

\end{document}